\pdfoutput=1
\documentclass[12pt]{amsart}

\usepackage[numbers]{natbib}
\usepackage[utf8]{inputenc}
\usepackage[margin=1in]{geometry}
\usepackage[dvipsnames, table]{xcolor}
\usepackage{color}
\usepackage{tikz}
\usetikzlibrary{backgrounds, positioning}
\usepackage{tikz-cd}
\usetikzlibrary{decorations.pathmorphing}
\usepackage{amssymb, amsmath, mathtools, amsfonts, amsthm}
\usepackage{graphicx}
\graphicspath{ {./images/} }
\usepackage{caption}
\usepackage{latexsym}
\usepackage[inline]{enumitem}
\usepackage{hyperref}
\hypersetup{colorlinks=true, urlcolor=caribbeangreen, citecolor=gray, linkcolor=caribbeangreen}
\usepackage{url}
\usepackage{indentfirst}
\usepackage{listings}
\usepackage{upquote}
\usepackage[mathscr]{eucal}
\usepackage{manfnt}
\usepackage{longtable}
\usepackage{colonequals}
\usepackage{nicefrac}
\usepackage{caption}
\usepackage{subcaption}
\usepackage{float}
\usepackage[capitalize]{cleveref}
\usepackage{stackengine}
\usepackage{breakcites}
\usepackage{stmaryrd}
\usepackage{bm}
\usepackage{soul}
\usepackage{comment}
\usepackage{booktabs}

\newtheorem{theorem}{Theorem}[section]
\newtheorem{lemma}[theorem]{Lemma}

\newtheorem{proposition}[theorem]{Proposition}

\theoremstyle{definition}
\newtheorem{definition}[theorem]{Definition}

\newtheorem{example}[theorem]{Example}
\newtheorem{remark}[theorem]{Remark}

\definecolor{santicolor}{RGB}{0, 159, 183}
\definecolor{nivencolor}{rgb}{0.9, 0.17, 0.31}
\definecolor{sameeracolor}{rgb}{0.57, 0.36, 0.51}
\definecolor{primitive}{HTML}{1A48FF}
\definecolor{defecttwo}{HTML}{F52773}

\definecolor{headercolor}{rgb}{0.83, 0.83, 0.83}
\definecolor{caribbeangreen}{rgb}{0.0, 0.8, 0.6}
\definecolor{acolor}{rgb}{0.01, 0.75, 0.24}
\definecolor{bcolor}{RGB}{0, 159, 183}
\definecolor{ccolor}{RGB}{254, 74, 73}
\definecolor{linkcolor}{rgb}{0.0, 0.0, 0.61}

\definecolor{magmabuiltin}{rgb}{0.10,0.20,0.45}
\definecolor{magmacomment}{gray}{0.45}
\definecolor{magmastring}{rgb}{0.10,0.40,0.25}
\definecolor{magmarule}{gray}{0.72}

\lstdefinelanguage{Magma}{%
  sensitive=true,
  morecomment=[l]{//},
  morecomment=[s]{/*}{*/},
  morestring=[b]",
  morekeywords=[1]{%
    adj,and,assert,assert2,assert3,assigned,break,by,case,cat,catch,clear,
    cmpeq,cmpne,continue,declare,default,delete,diff,div,do,elif,else,end,eq,
    error,eval,exists,exit,false,for,forall,forward,fprintf,freeze,function,
    ge,gt,if,iload,import,in,intrinsic,is,join,le,load,local,lt,meet,mod,ne,
    not,notadj,notin,notsubset,or,print,printf,procedure,quit,random,read,
    readi,repeat,require,requirege,requirerange,restore,return,save,sdiff,
    select,subset,then,time,to,true,try,until,vprint,vprintf,vtime,when,
    where,while,xor,
    car,cop,elt,ext,hom,ideal,iso,lideal,map,ncl,pmap,quo,rec,recformat,
    rideal,sub},
  morekeywords=[2]{%
    Integers,Rationals,Reals,Complexes,IntegerRing,RationalField,RingOfIntegers,
    FiniteField,GF,NumberField,CyclotomicField,PolynomialRing,FunctionField,
    FieldOfFractions,AffineSpace,ProjectiveSpace,Scheme,Curve,EllipticCurve,
    HyperellipticCurve,Jacobian,GenusOneModel,cInvariants,bInvariants,
    Invariants,MinimalModel,WeierstrassModel,jInvariant,Discriminant,Conductor,
    Rank,Regulator,TorsionSubgroup,MordellWeilGroup,Generators,RationalPoints,
    Points,Genus,Degree,Dimension,BaseChange,BaseRing,Parent,Factorization,
    Factorisation,GCD,LCM,Valuation,Sqrt,Roots,Coefficients,Evaluate,Numerator,
    Denominator,Matrix,Vector,Eltseq,Domain,Codomain,Kernel,Image,Order,Index,
    IsSquare,IsPrime,IsIrreducible,IsIsomorphic,IsSingular,Automorphisms,
    AutomorphismGroup,DefiningEquations,Ideal,Set,Sequence,SetToSequence,
    Sprintf,Random,Abs,Max,Min,Support,Divisor},
}

\lstdefinestyle{magmastyle}{%
  language=Magma,
  basicstyle=\ttfamily\footnotesize,
  keywordstyle=[1]\bfseries,
  keywordstyle=[2]\color{magmabuiltin},
  commentstyle=\itshape\color{magmacomment},
  stringstyle=\color{magmastring},
  columns=fullflexible,
  keepspaces=true,
  upquote=true,
  showspaces=false,
  showstringspaces=false,
  showtabs=false,
  tabsize=2,
  breaklines=true,
  breakatwhitespace=true,
  postbreak=\mbox{\textcolor{magmacomment}{$\hookrightarrow$}\space},
  frame=tb,
  framerule=0.4pt,
  rulecolor=\color{magmarule},
  framesep=6pt,
  xleftmargin=\parindent,
  aboveskip=1.1em,
  belowskip=1.1em,
  numbers=none,
  captionpos=b,
  escapeinside={(*}{*)},
}

\lstnewenvironment{magma}[1][]{\lstset{style=magmastyle,#1}}{}

\newcommand{\R}{\mathbb{R}}
\newcommand{\Z}{\mathbb{Z}}

\newcommand{\Q}{\mathbb{Q}}
\newcommand{\QQ}{\Q}

\newcommand{\Qbar}{{\bar{\Q}}}

\newcommand{\OO}{\mathscr{O}}

\newcommand{\A}{\mathbb{A}}
\newcommand{\Gm}{\mathbb{G}_\text{m}}
\newcommand{\Gw}[1][]{\mathbb{G}_{\bfw  #1}}
\newcommand{\kk}{{\mathbf k}}
\renewcommand{\P}{\mathbb{P}}
\newcommand{\Pone}{\P^1}

\newcommand{\mcX}{\mathcal{X}}

\newcommand{\mcU}{\mathcal{U}}

\newcommand{\mcF}{\mathcal{F}}

\newcommand{\mcL}{\mathcal{L}}

\newcommand\HH{\mathrm{H}}

\newcommand{\x}{\mathsf{x}}
\newcommand{\y}{\mathsf{y}}
\newcommand{\z}{\mathsf{z}}
\newcommand{\sfu}{\mathsf{u}}

\newcommand{\pr}{\mathrm{pr}}
\renewcommand{\t}{\mathsf{t}}

\newcommand{\sfK}{\mathsf{K}}

\newcommand{\chiabc}{\tfrac1a+\tfrac1b+\tfrac1c-1}

\renewcommand{\leq}{\leqslant}
\renewcommand{\geq}{\geqslant}

\renewcommand{\emptyset}{\varnothing}

\newcommand{\paren}[1]{\left( #1 \right)}
\newcommand{\cdef}[1]{{\color{linkcolor}\textsf{#1}}}
\newcommand{\brk}[1]{\left\{ #1 \right\}}

\newcommand{\bfe}{\mathbf{e}}

\newcommand{\bfw}{\mathbf{w}}
\newcommand{\gfe}{A\x^{a} + B\y^{b} + C\z^{c}}
\newcommand{\ggfe}{c_1\x_1^{e_1} + c_2\x_2^{e_2} + c_3\x_3^{e_3}}
\newcommand{\eee}{(e_1, e_2, e_3)}
\newcommand{\chie}{\tfrac{1}{e_1} + \tfrac{1}{e_2} + \tfrac{1}{e_3} - 1}

\newcommand{\Autsch}{\mathbf{Aut}}

\newcommand{\Reg}{\operatorname{Reg}}

\newcommand{\bad}{\mathsf{B}}

\newcommand{\Magma}{\texttt{Magma}}

\DeclareMathOperator{\Spec}{Spec}
\DeclareMathOperator{\Proj}{Proj}
\DeclareMathOperator{\Gal}{Gal}
\DeclareMathOperator{\Stab}{Stab}

\DeclareMathOperator{\Aut}{Aut}
\DeclareMathOperator{\Jac}{Jac}

\DeclareMathOperator{\Pic}{Pic}

\DeclareMathOperator{\lcm}{lcm}
\DeclareMathOperator{\Ht}{Ht}
\DeclareMathOperator{\logHt}{ht}

\DeclareMathOperator{\covol}{covol}

\DeclareMathOperator{\tors}{tors}
\DeclareMathOperator{\defect}{defect}

\usepackage[textsize=tiny]{todonotes}
\usepackage{xargs, xpatch}
\usepackage{marginnote}

\usetikzlibrary{calc,intersections}

\definecolor{tilegray}{gray}{0.42}

\tikzset{
  tileface/.style = {draw=black, line width=0.3pt},
  tiledark/.style = {tileface, fill=tilegray},
  tilelite/.style = {tileface, fill=white},
}

\newsavebox\tileDthree
\begin{lrbox}{\tileDthree}
\begin{tikzpicture}[scale=1]
  \def\R{1.15}
  \coordinate (O) at (0,0);
  \foreach \k in {0,...,5} {\coordinate (V\k) at ({60*\k}:\R);}

  \foreach \k [evaluate=\k as \kk using {int(mod(\k+1,6))}] in {0,2,4}
    \filldraw[tiledark] (O) -- (V\k) -- (V\kk) -- cycle;
  \foreach \k [evaluate=\k as \kk using {int(mod(\k+1,6))}] in {1,3,5}
    \filldraw[tilelite] (O) -- (V\k) -- (V\kk) -- cycle;

\end{tikzpicture}
\end{lrbox}

\newsavebox\tileDsix
\begin{lrbox}{\tileDsix}
\begin{tikzpicture}[scale=1]
  \def\R{1.15}
  \coordinate (O) at (0,0);
  \foreach \k in {0,...,5} {\coordinate (V\k) at ({60*\k}:\R);}
  \foreach \k in {0,...,5} {\coordinate (M\k) at ({60*\k+30}:{\R*sqrt(3)/2});}

  \foreach \k [evaluate=\k as \kk using {int(mod(\k+1,6))}] in {0,...,5}{
    \filldraw[tiledark] (O) -- (V\k) -- (M\k)  -- cycle;
    \filldraw[tilelite] (O) -- (M\k) -- (V\kk) -- cycle;
  }

\end{tikzpicture}
\end{lrbox}

\newsavebox\tileDfour
\begin{lrbox}{\tileDfour}
\begin{tikzpicture}[scale=1]
  \def\S{1.0}
  \coordinate (O) at (0,0);
  \coordinate (W0) at (\S,\S);
  \coordinate (W1) at (-\S,\S);
  \coordinate (W2) at (-\S,-\S);
  \coordinate (W3) at (\S,-\S);

  \foreach \k [evaluate=\k as \kk using {int(mod(\k+1,4))}] in {0,2}
    \filldraw[tiledark] (O) -- (W\k) -- (W\kk) -- cycle;
  \foreach \k [evaluate=\k as \kk using {int(mod(\k+1,4))}] in {1,3}
    \filldraw[tilelite] (O) -- (W\k) -- (W\kk) -- cycle;

\end{tikzpicture}
\end{lrbox}

\newcommand{\tessfigure}{%
\begin{figure}[htbp]
  \centering
  \begin{tabular}{@{}c@{\hspace{2.4em}}c@{\hspace{2.4em}}c@{}}
    \usebox\tileDthree & \usebox\tileDsix & \usebox\tileDfour \\[0.8ex]
    $\{3,3,3\}$ & $\{3,2,6\}$ & $\{4,2,4\}$ \\
  \end{tabular}
  \caption{The angles of each triangle are $\pi/e_1, \pi/e_2, \pi/e_3$.}
  \label{fig:euclidean-tiles}
\end{figure}}

\title[Counting Euclidean solutions]{Counting solutions to Euclidean Fermat equations}

\author{Niven T. Achenjang}
\address{Department of Mathematics, Harvard University, Cambridge, MA 02138, USA}
\email{achenjang@math.harvard.edu }
\urladdr{\url{https://people.math.harvard.edu/~achenjang/}}

\author{Santiago Arango-Pi{\~n}eros}
\address{Department of Mathematics, University of Massachusetts Amherst, Amherst, MA 01003, USA}
\email{santiago.arango.pineros@gmail.com}
\urladdr{\url{https://sarangop1728.github.io/}}

\author{Sameera Vemulapalli}
\address{Department of Mathematics, Harvard University, Cambridge, MA 02138, USA}
\email{vemulapalli@math.harvard.edu}
\urladdr{\url{https://sites.google.com/view/sameeravemulapalli/}}

\begin{document}

\begin{abstract}
  A solution $(x,y,z) \in \mathbb{Z}^3-\brk{(0,0,0)}$ to a Fermat equation
\begin{equation}
    \label{eq:abstract-gfe}
    \gfe = 0
\end{equation}
is called \emph{primitive} if $\gcd(x,y,z) = 1$. In the \emph{Euclidean} regime
(that is, when the Euler characteristic
$\chi = \tfrac{1}{a} + \tfrac{1}{b} + \tfrac{1}{c} - 1$ is zero),
\Cref{eq:abstract-gfe} defines a smooth genus one curve in a weighted projective
plane. We study these models and provide an asymptotic formula for the number
of primitive integral solutions with bounded height.
\end{abstract}

\maketitle
\setcounter{tocdepth}{1}
\tableofcontents
\vspace{-11mm}

\section{Introduction}\label{sec:introduction}
Fix integers $a,b,c\geq 2$ and $A,B,C$ with $ABC \neq 0$. An integer solution $(x,y,z)$ to the \cdef{Fermat equation}
\begin{equation}
\label{eq:gfe}
 F \colonequals \gfe \in \Z[\x,\y,\z], \quad F = 0
\end{equation}
is called \cdef{primitive} if $\gcd(x,y,z) = 1$. The \cdef{multiplicative Weil
  height of $\Pone(\Q)$} is the function $\Ht\colon \Pone(\Q) \to \Z_{> 0}$ given by
$\Ht(Q) = \max\brk{|s|,|t|},$ where $s,t \in \Z$ are coprime integers such that
$(s : t) = Q$. For a subset $\Omega \subseteq \Pone(\Q)$ and any positive real number
$H$, define
\begin{align}
  N(\Omega; H) \colonequals \#\brk{Q \in \Omega : \Ht(Q) \leq H}.
\end{align}

The purpose of this article is to prove the following: 

\begin{theorem}
  \label{thm:main}
  Suppose that \cref{eq:gfe} is \cdef{Euclidean}, that is: $\chiabc = 0$. Let
  $\Omega_F$ be the subset of rational points $(-Ax^a:Cz^c) \in \Pone(\Q)$ arising
  from primitive integral solutions $(x,y,z)$ to \cref{eq:gfe}, and suppose
  that $\Omega_F \neq \emptyset$. Then there is an integer $r(F) \geq 0$
  (\Cref{remark:leading-constant}) such that the following hold.
  \begin{enumerate}[leftmargin=*, label=(\alph*)]
  \item \label{it:main-rank-zero} If $r(F) = 0$, then $\Omega_F$ is finite, and
    $N(\Omega_F;H) = \#\Omega_F$ for all $H$ sufficiently large.
  \item \label{it:main-positive-rank} If $r(F) > 0$, then there is a constant $\kappa(F) > 0$ such that
    \begin{equation*}
      N(\Omega_F;H) = \kappa(F)\cdot \paren{\log H}^{r(F)/2}
      + O\paren{\paren{\log H}^{(r(F)-1)/2}}, \qquad H \to \infty.
    \end{equation*}
  \end{enumerate}
\end{theorem}

\begin{remark}
  \label{remark:leading-constant}
  Associated to \cref{eq:gfe} is a scheme $\mcX$ defined by the same
  equation in an appropriate weighted projective space over $\Z$; see
  \Cref{sec:generalities}. The generic fiber $X \coloneqq \mcX_{\QQ}$ is a smooth genus $1$
  curve, and we let $E\coloneqq\Jac(X)$ denote its Jacobian. 
  When $\Omega_F \neq \emptyset$, the curve $X$ has a rational point and 
  so is (non-canonically) isomorphic to $E$. The
  quantity $r(F)$ in \cref{thm:main} is the \cdef{Mordell--Weil rank} of
  $E(\Q)$. Furthermore, when $r(F) >0$, the constant is given by
  \begin{equation}
    \label{eq:kappa-F}
    \kappa(F) = \varrho_F \cdot
    \frac{\#E(\Q)_{\tors}}
    {\#\Aut(\phi_{\Q})\cdot\Gamma\paren{1 + \tfrac{r(F)}{2}}\cdot \sqrt{\Reg(E/\Q)}}
    \cdot \paren{\dfrac{2\pi\lcm(a,b,c)}{abc}}^{r(F)/2},
  \end{equation}
  where $\Gamma(s)$ is the \cdef{Gamma function}, $E(\Q)_{\tors}$ is the
  \cdef{torsion subgroup} of $E(\Q)$, $\Reg(E/\Q)$ is the \cdef{regulator} with respect to the line bundle $\OO_E(2O)$, where $O$ is the origin of $E$,
  $\Aut(\phi_\Q)$ is the \cdef{group of $\Q$-automorphisms} of the map
  $\phi\colon X \to \Pone_{\Q}$, taking $(x,y,z) \mapsto (-Ax^a:Cz^c)$
  (\Cref{lemma:aut-rational}), and $\varrho_F \in (0,1]$ is the proportion of rational
  points of $X$ that come from \emph{primitive} integral solutions
  (\cref{eq:rho}).
\end{remark}

The proof of \Cref{thm:main} has two main ingredients. The analytic component
is a classical theorem of Néron \cite{Neron-65}, which computes the asymptotic
distribution of rational points with bounded canonical height on an abelian
variety. This is already sufficient to prove our theorem in the case of
Selmer-type cubics $A\x^3+ B\y^3 + C\z^3 = 0$ (since every rational point comes
from two primitive integral solutions).

The second ingredient is a sieve for the \emph{primitive
    locus} (see \Cref{thm:primitive-subgroup}): not all points of $X(\Q)$ arise
  from primitive solutions. We show that the primitive locus is a finite union
  of cosets of a finite index subgroup of the Mordell--Weil group of (any
  elliptic curve structure on) $X$, and then apply N\'eron's theorem. This
  \emph{primitivity sieve} is the crux of the general argument and our main
  original contribution. We explain it via an example.

\begin{example}
  \label{ex:intro-defect} Consider the Fermat equation
  \begin{equation}\label{eq:gfe-example}
    F \colonequals \x^3 - 2\y^2 + \z^6,\quad F=0.
  \end{equation}
  Let $X \subset \P(2,3,1)_\Q$ be the associated curve. It is a nice curve
  of genus one. For every $P \in X(\Q)$ we can find a representative
  $(x,y,z) \in \Z^3-\{(0,0,0)\}$ that is minimal with respect to the weighted
  scaling
  $(x,y,z)\cdot \lambda \colonequals (\lambda^2 x, \lambda^3 y,\lambda z)$ for some
  $\lambda \in \Q^{\times}$. That is, $P = (x:y:z)$ with $x,y,z \in \Z$ such that no prime
  $p$ satisfies
  $p^2 \mid x, \, p^3 \mid y, \,p \mid z$ simultaneously. Moreover, any two such
  representatives $(x_1,y_1,z_1)$ and $(x_2,y_2,z_2)$ may differ only by the
  action of $-1$:
  $(x_2,y_2,z_2) = ((-1)^2x_1,(-1)^3y_1,-z_1) = (x_1,-y_1,-z_1)$. In
  particular, we may define the \cdef{primitivity defect} of a rational point
  $P \in X(\Q)$ to be $\gcd(x,y,z)$ for any minimal representative $(x,y,z)$. It
  is straightforward to show that the subset $X(\Q)_1 \subseteq X(\Q)$ of points
  ``coming from primitive integral solutions'' to \cref{eq:gfe-example} consists
  precisely of the points of defect one.

  Choose $O = (2:2:0) \in X(\Q)$ as the origin, and let $E = (X,O)$ be the
  resulting elliptic curve. The point $P = (1:1:1)$ is primitive, but already
  $2P$ has minimal representative $(-14,26,4)$, of defect $2$;
  \Cref{tab:intro-defect} shows the defects of the first few multiples of $P$
  alternating between $1$ and $2$. The pattern becomes transparent in a basis:
  $E(\Q) \cong \Z \oplus \Z/2\Z$, generated by $Q = (2:6:2)$ of infinite order and
  the $2$-torsion point $T = (-1:0:1)$, with $P = T - Q$. For all $|m| \leq 40$,
  the point $mQ$ has defect $2$ and the point $mQ + T$ has defect $1$
  (see \Cref{fig:intro-torsion}), suggesting that
  \[
    X(\Q)_1 = T + \Lambda, \qquad \Lambda = \langle Q \rangle.
  \]
  That the primitive locus is a finite union of cosets of a finite index
  subgroup is no accident: this is the content of
  \Cref{thm:primitive-subgroup-intro}. One can prove that the density of
  \Cref{def:rho} is $\varrho_F = 1/2$, but, to keep this example brief, we omit the
  details here.
\begin{table}[ht]
  \centering
  \caption{The first nine multiples of $P = (1:1:1)$ on
    $X\colon \x^3 - 2\y^2 + \z^6 = 0$.}
  \footnotesize
  \setlength{\tabcolsep}{4pt}
  \rowcolors{2}{white}{gray!15}
  \begin{tabular}{clc}
    \toprule
    $m$ & Minimal representative $(x,y,z)$ for $mP = (x:y:z)$ & $\gcd(x,y,z)$ \\
    \midrule
    $1$ & $(1 , 1 , 1)$ & $1$ \\
    $2$ & $(-14 , 26 , 4)$ & $2$ \\
    $3$ & $(17 , -2387 , 15)$ & $1$ \\
    $4$ & $(62146 , -10983646 , 104)$ & $2$ \\
    $5$ & $(10374961 , 23894024939 , 1711)$ & $1$ \\
    $6$ & $(-3098082958 , 2073702872349862 , 143220)$ & $2$ \\
    $7$ & $(-34197096105119 , -184226381697790488601 , 6899489)$ & $1$ \\
    $8$ & $(1785866978837560066 , -1989644569090679752355375234 , 1142299184)$ & $2$ \\
    $9$ & $(94615651228339230126497 , 20579268644229142617017417250997993 , 34859070945)$ & $1$ \\
    \bottomrule
  \end{tabular}
  \label{tab:intro-defect}
\end{table}

\begin{figure}[ht]
  \centering
  \includegraphics[width=0.7\textwidth]{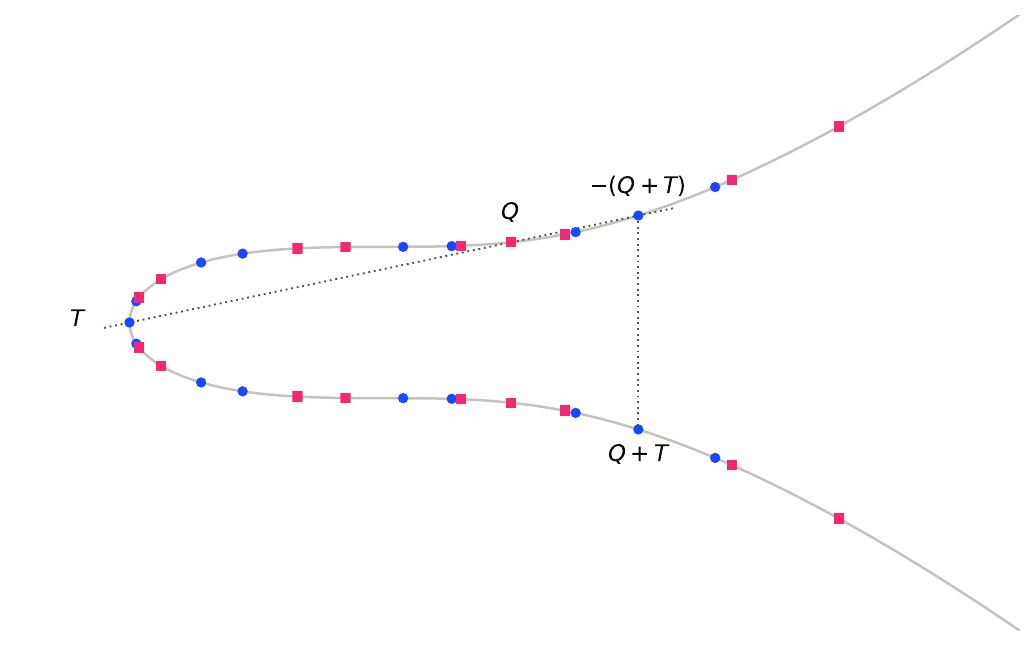}
  \caption{Rational points on the curve
    $X\colon \x^3-2\y^2+\z^6=0$, drawn in the coordinates of
    $E\colon y^2 = x^3 + 8$, colored by primitivity defect:
    ${\color{primitive} X(\Q)_1}$ (blue circles) and ${\color{defecttwo} X(\Q)_2}$
    (pink squares). Translation by $T$ swaps the two
    colors.}
  \label{fig:intro-torsion}
\end{figure}
\end{example}

Our sieving result states that for any Euclidean Fermat equation with a non-trivial integral solution, the primitive locus $X(\Q)_1$ of the corresponding weighted projective curve $X$ is a finite union of orbits of a finite index subgroup of $E(\QQ)$.

\begin{theorem}[The Primitive Subgroup]
  \label{thm:primitive-subgroup-intro}
  Suppose that \cref{eq:gfe} is Euclidean and that $X(\Q) \neq \emptyset$.
  Let $E\colonequals\Jac(X)$. Then the subgroup
  \[
    \Lambda \colonequals \brk{P \in E(\Q) : P + X(\Q)_1 = X(\Q)_1}
  \]
  has finite index in $E(\Q)$, and $X(\Q)_1$ is a finite disjoint union of
  orbits of $\Lambda$ acting on $X(\Q)$. Moreover, $\Lambda$ is the largest subgroup of $E(\Q)$
  with this property.
  We call it the \cdef{primitive subgroup}.
\end{theorem}

\subsection{Previous work}
\label{sec:previous-work} As Darmon and Granville observed in~\cite[Section~6]{Darmon&Granville95},
the set of primitive integral solutions to a Euclidean Fermat equation is
either finite or infinite, as determined by the ranks
of a finite list of twists of certain CM elliptic curves; see 
\cite{Beukers98, BCDY15,BennetMihlailescuSiksek16,WG24} for existing results on Euclidean signatures contextualized in the broader generalized Fermat literature. Most of the research on this subject is focused on the rational points of these twist families of genus one curves rather than the primitive integral solutions of the underlying generalized Fermat equation.

For signature $(3,3,3)$ the classical reference is Selmer's study of the cubics
$A\x^3+B\y^3+C\z^3=0$ \cite{Selmer-51,Selmer-54}. In this foundational
computational work, the descent that inspired the definition of the
\emph{Selmer group} was carried out, and the first cubic counterexample to
the Hasse principle was found:
\[3\x^3+4\y^3+5\z^3=0.\]
See also Cassels \cite{Cassels66} and Mordell \cite{Mordell69}. The diagonal
subfamily $\x^3+\y^3=n\z^3$, the ``sum of two rational cubes'' problem, has
its own history, running from Sylvester \cite{Sylvester79} through Satg\'e's
Heegner-type constructions \cite{Satge86,Satge87} and Dasgupta--Voight
\cite{DV18} to the density results of Alp\"oge--Bhargava--Shnidman
\cite{ABS22} and Koymans--Smith \cite{KoymansSmith:SumsOfCubes}.

For signature $(4,2,4)$, an even older counterexample to the Hasse principle for
genus one curves was found by Lind \cite{Lind40} and Reichardt
\cite{Reichardt42}:
\[2\y^2=\x^4-17\z^4.\]
See \cite[Thm.~X.6.5]{Silverman09} for a modern proof, and
\cite{Lemmermeyer99,AL11} for an account of earlier work on equations of
this shape going back to P\'epin and Mordell.

For signature $(3,2,6)$, the curve $X : F = 0 \subset \P(2,3,1)_{\Q}$ has the
rational point at infinity $\z=0$, so it is an elliptic curve (a sextic twist
of $\y^2=\x^3+1$) and the dichotomy above is decided by its rank alone, provided that a primitive integral solution exists. These
are the Mordell curves $\y^2=\x^3+k$, whose rational points were studied by
Cassels \cite{Cassels50} and Selmer \cite{Selmer-56} by descent; for the
history of the problem and extensive computations see Mordell
\cite[Ch.~26]{Mordell69} and Gebel--Peth\H{o}--Zimmer \cite{GPZ98}. More recently,
the arithmetic statistics of this family ($3$-isogeny Selmer groups, ranks)
have been studied by Bhargava--Elkies--Shnidman \cite{BES20}.

Our focus in this paper is different from all of these directions. We
restrict our attention to those Fermat equations admitting a primitive
integral solution, and study the asymptotic size of the set of solutions,
following the narrative of \cite{ArangoPineros-CountingSphericalSolutions},
which treats the spherical analogue of this problem.

\subsection*{Acknowledgments}
\label{sec:acknowledgements} 

The first author was supported by the NSF under grant number DMS-2503296.
The third author was supported by the NSF under grant number DMS-2303211. SAP
thanks Louis Gaudet for stimulating questions and conversations related to
this project.

\subsection*{AI disclosure}
\label{sec:ai-disclosure}

The mathematical ideas here arose from human collaboration. We used AI tools
to create examples and Figure 1, to produce a referee report on an older draft,
and to carry out a literature review. We wrote this paper using our own words.

\section{Generalities}
\label{sec:generalities}

\subsection{Setup} We establish some notation associated to \Cref{eq:gfe}. For
notational convenience, we rewrite the polynomial $F$ of \Cref{eq:gfe} as
$F \colonequals \ggfe \in \Z[\x_1,\x_2,\x_3]$.

\begin{itemize}[leftmargin=*]
\item The symbols $i,j$ will denote indices, ranging in $\{1,2,3\}$.
\item The coefficients $c_i \in \Z$ satisfy $c_1c_2c_3\neq0$.
\item The exponents $e_i \in \Z_{\geq 2}$. The triple $\bfe \colonequals \eee$ is
  called a \cdef{signature} of $F$.
\item The quantity $\chi \colonequals \chie$ is called the \cdef{Euler
    characteristic} of the signature $\bfe$.
\item The \cdef{weight vector} corresponding to the signature $\bfe$ is
  \[
    \bfw \colonequals (w_1, w_2, w_3) \colonequals \left(\tfrac{\ell}{e_1}, \tfrac{\ell}{e_2},
    \tfrac{\ell}{e_3} \right)\in \Z_{\geq 1}^{3},
  \]
  where $\ell \colonequals \lcm\eee$. From the identity
  \begin{equation}
    \label{eq:triple-lcm}
    e_1e_2e_3 = \ell\cdot\gcd(e_2e_3,e_1e_3,e_1e_2),
  \end{equation} we see that $\gcd(\bfw) = 1$.
\item We say that the weight vector $\bfw$ is \cdef{well-formed} if the
  weights are pairwise relatively prime: $\gcd(w_i,w_j) = 1$ for all $i \neq j$.
\item We consider the polynomial ring $\Z[\x_1,\x_2,\x_3]$ with grading given
  by $\deg(\x_i) = w_i$.
\item $\P(\bfw) \colonequals \Proj \Z[\x_1,\x_2,\x_3]$ is the weighted projective space
  associated to the graded polynomial ring defined above.
\item $\mcX \colonequals \Proj \Z[\x_1,\x_2,\x_3]/(F) \subset \P(\bfw)$ is the weighted
  projective scheme associated to \Cref{eq:gfe}.
\item $X \colonequals \mcX_\Q$ is the generic fiber of $\mcX$.
\item $\mcU \colonequals \{F=0\}\subset\A^3_\Z-\mathbf0$ is the punctured cone above $\mcX$.
\item $\Gm \colonequals \Spec \Z[\t,\t^{-1}]$ is the multiplicative group.
\item $\mu_n \colonequals \Spec \Z[\t, \t^{-1}]/(\t^n-1)$ is the $n$-torsion
  subgroup scheme of $\Gm$.
\item $\Gw$ is the scheme-theoretic image of the injective morphism of
  group schemes
  \[\Gm \to \Gm^3, \quad t \mapsto (t^{w_1}, t^{w_2}, t^{w_{3}}).\]
  In particular, $\Gw \cong \Gm$.
\end{itemize}

\subsection{The weighted multiplicative action}
\label{sec:weighted-multiplicative-action}
The group $\Gw$ acts on $\mcU$ via 
\begin{equation*}
  (x_1,x_2,x_3)\cdot t \colonequals (t^{w_1}x_1, t^{w_2}x_2, t^{w_3}x_3),
\end{equation*}
and we understand the stabilizers of this action.

\begin{lemma}
  \label[lemma]{lemma:stabilizers}
  Let $K$ be a field and let $P = (x_1, x_2, x_3) \in \mcU(K)$. Then, we compute the
  stabilizer subgroup of $P$ in $\Gw[,K]$ to be 
  \[
    \Stab_{\Gw[,K]}(P) = \bigcap_{x_i \, \neq \, 0} \mu_{w_i, K}.
  \]
  In particular:
  \begin{enumerate}[label=(\alph*), leftmargin=*]
  \item \label{it:one-nonzero} If exactly one coordinate $x_i$ is nonzero, then
    $c_i = 0$ in $K$, and the stabilizer is $\mu_{w_i, K}$. 
  \item \label{it:two-nonzero} If $\bfw$ is well-formed and at least two
    coordinates of $P$ are nonzero, the stabilizer is trivial.
  \end{enumerate}
\end{lemma}
\begin{proof}
  A point $t$ of $\Gw$ fixes $P$ if and only if $t^{w_i} = 1$ for every $i$
  with $x_i \neq 0$; this is the displayed intersection. For \ref{it:one-nonzero},
  the defining equation of $\mcU$ collapses to $c_i x_i^{e_i} = 0$ with
  $x_i \neq 0$, forcing $c_i = 0$. For \ref{it:two-nonzero}, if
  $x_i, x_j \neq 0$ with $i \neq j$ then
  $\Stab_{\Gw}(P) = \mu_{w_i} \cap \mu_{w_j} = \mu_{\gcd(w_i,w_j)}$, which is trivial
  because $\bfw$ is well-formed.
\end{proof}

\begin{definition}
    \label{def:bad-locus}
    The \cdef{bad locus} of \cref{eq:gfe} is the integer 
    \begin{equation}
  \label{eq:bad}
  \bad_F \colonequals \prod_{\substack{i \in \{1,2,3\} \\ w_i > 1}} c_i,
\end{equation}
an empty product being $1$. We say that a prime is \cdef{good} if it does not divide $\bad_F$. We will work with the ring
\begin{equation}
  \label{eq:R}
  R \colonequals R_F \colonequals \Z[1/\bad_F].
\end{equation}
\end{definition}

\begin{lemma}
  \label[lemma]{lemma:U->X-Gm-torsor}
  Suppose the weight vector $\bfw$ is well-formed.
  Then, the morphism $\pr\colon \mcU_R \to \mcX_R$ given by $(x_1,x_2,x_3) \mapsto (x_1:x_2:x_3)$ is a
  $\Gw[,R]$-torsor. In particular, the induced map on $R$-points
  $\mcU(R) \to \mcX(R)$ is $|R^\times|$-to-one and surjective.
\end{lemma}
\begin{proof}
  Let $I\coloneqq\{i\in\{1,2,3\}\colon w_i=1\}$. Then, by definition of
  $\bad_F$, the following is an open cover of $\mcX_R$:
    \[\mcX_R=D_+(\x_1\x_2)\cup D_+(\x_1\x_3)\cup D_+(\x_2\x_3)\cup\bigcup_{i\in I}D_+(\x_i).\]
    Hence, it suffices to check that $\pr\colon\mcU_R\to \mcX_R$ is a $\Gw[,R]$-torsor
    over each set in the above open cover.

    We first argue over $D_+(\x_i)$, where $i\in I$. Let $S\coloneqq R[\x_1,\x_2,\x_3]/(F)$, and let $S_{(\x_i)}$ denote the degree $0$ part of $S[1/\x_i]$. Over $D_+(\x_i)\cong\Spec S_{(\x_i)}$, the map $\pr\colon\mcU\to \mcX$ becomes the map $\Spec S[1/\x_i]\to\Spec S_{(\x_i)}$ induced by the natural inclusion $S_{(\x_i)}\hookrightarrow S[1/\x_i]$. We need to show this is a $\Gw$-torsor, i.e., for any $R$-algebra $A$ and $\phi_0\colon S_{(\x_i)}\to A$, $\Gw(A)$ acts simply transitively on the set of $R$-algebra maps $\phi\colon S[1/\x_i]\to A$ which restrict to $\phi_0$. Indeed, for two such ring maps $\phi,\phi'\colon S[1/\x_i]\to A$, set $\lambda\coloneqq\phi(\x_i)/\phi'(\x_i)\in A^\times$ and observe that, for $j\in\{1,2,3\}$,
    \[\phi(\x_j)=\phi(\x_i)^{w_j}\phi\paren{\frac{\x_j}{\x_i^{w_j}}}=\phi(\x_i)^{w_j}\phi_0\paren{\frac{\x_j}{\x_i^{w_j}}}=\phi(\x_i)^{w_j}\phi'\paren{\frac{\x_j}{\x_i^{w_j}}}=\lambda^{w_j}\phi'(\x_j),\]
    so $(\lambda^{w_1},\lambda^{w_2},\lambda^{w_3})\in\Gw(A)$ is the unique element sending $\phi'$ to $\phi$, finishing the argument for this open.

    We finish by arguing over $D_+(\x_1\x_2)$; the arguments over the remaining opens are analogous. Let $S$ be as before, and let $S_{(\x_1\x_2)}$ denote the degree $0$ part of $S[1/\x_1\x_2]$. We need to show that, for any $R$-algebra $A$, given $\phi,\phi'\colon S[1/\x_1\x_2]\to A$ which agree on $S_{(\x_1\x_2)}$, there is a unique element of $\Gw(A)$ sending $\phi'$ to $\phi$. Fixing such $\phi,\phi'$, first observe that
    \[\frac{\phi(\x_2)^{w_1}}{\phi(\x_1)^{w_2}}=\phi\paren{\frac{\x_2^{w_1}}{\x_1^{w_2}}}=\phi'\paren{\frac{\x_2^{w_1}}{\x_1^{w_2}}}=\frac{\phi'(\x_2)^{w_1}}{\phi'(\x_1)^{w_2}}\in A.\]
    Set $\delta_1\coloneqq\phi(\x_1)/\phi'(\x_1)\in A^\times$, so the above shows that $\phi(\x_2)^{w_1}=\phi'(\x_2)^{w_1}\cdot\delta_1^{w_2}$. Setting $\delta_2\coloneqq\phi(\x_2)/\phi'(\x_2)\in A^\times$, one similarly shows that 
    \begin{equation}\label{eqn:weighted-ratios}
        \phi(\x_3)^{w_1}=\phi'(\x_3)^{w_1}\cdot\delta_1^{w_3},\quad\phi(\x_3)^{w_2}=\phi'(\x_3)^{w_2}\cdot\delta_2^{w_3},\quad\text{and }\phi(\x_1)^{w_2}=\phi'(\x_1)^{w_2}\cdot\delta_2^{w_1}.
    \end{equation}
    The last identity says that $\delta_1^{w_2}=\delta_2^{w_1}$. Since $\gcd(w_1,w_2)=1$, we may choose $n,m\in\Z$ such that $mw_2-nw_1=-1$. Set
    \[\lambda\coloneqq\delta_1^{n}\delta_2^{-m}\in A^\times.\]
    Then $\lambda^{w_1}=\delta_1^{nw_1}\delta_2^{-mw_1}=\delta_1^{nw_1}\delta_1^{-mw_2}=\delta_1$ and $\lambda^{w_2}=\delta_2^{nw_1}\delta_2^{-mw_2}=\delta_2$. Moreover, $\lambda$ is the only element of $A^\times$ with these two properties, since any such element $\mu$ satisfies $\mu=\mu^{nw_1-mw_2}=\delta_1^n\delta_2^{-m}$. To finish the proof, we must show that $\phi(\x_3)=\phi'(\x_3)\cdot\lambda^{w_3}$, for then $(\lambda^{w_1},\lambda^{w_2},\lambda^{w_3})\in\Gw(A)$ will be the unique element sending $\phi'$ to $\phi$. For this, we observe
    \[\frac{\phi(\x_3)\phi(\x_2)^{w_3 m}}{\phi(\x_1)^{w_3 n}}=\phi\paren{\frac{\x_3\x_2^{w_3 m}}{\x_1^{w_3 n}}}=\phi'\paren{\frac{\x_3\x_2^{w_3 m}}{\x_1^{w_3 n}}}=\frac{\phi'(\x_3)\phi'(\x_2)^{w_3 m}}{\phi'(\x_1)^{w_3 n}},\]
    so
    \[\phi(\x_3)=\phi'(\x_3)\cdot\paren{\frac{\phi'(\x_2)^m\phi(\x_1)^n}{\phi'(\x_1)^n\phi(\x_2)^m}}^{w_3}=\phi'(\x_3)\cdot\paren{\delta_1^n/\delta_2^m}^{w_3}=\phi'(\x_3)\cdot\lambda^{w_3}.\qedhere\]
\end{proof}

\begin{remark}
  In general, the projection map $\pr\colon \mcU \to \mcX$ is not a
  $\Gw$-torsor. For example, for the equation $\x_1^2+\x_2^3+\x_3^7=0$, we have $\bfw = (21,14,6)$, which is not well-formed. Since $\mcX \cong \Pone$, $\mcX(\Z) \cong X(\Q)$ is infinite. On the other hand,
  $\#\mcU(\Z) = 16$ by \cite[Theorem 1.1]{PoonenSchaeferStoll07}. This means
  that the map $\mcU(\Z) \to \mcX(\Z)$ is not surjective, and thus $\pr\colon \mcU \to
  \mcX$ cannot be a $\Gw$-torsor.
\end{remark}

\begin{lemma}
  \label[lemma]{lemma:genus}
  The scheme $X$ is a smooth, projective, and geometrically irreducible
  curve of genus
  \begin{equation}
    \label{eq:genus}
    g(X) = 1 - \tfrac12 \left(\gcd(e_2,e_3) + \gcd(e_1,e_3) + \gcd(e_1,e_2) - \gcd(e_2e_3,e_1e_3,e_1e_2) \right).
  \end{equation}
\end{lemma}
\begin{proof}
  Using \cite[Lemma 5.3.4]{Voight&Zureick-Brown22}, we conclude that $X$
  is smooth since it is the coarse space of the stack $[\mcU_{\Q}/\Gw[,\Q]]$.
  In general, weighted projective varieties are projective.

  By B\'ezout's theorem for weighted projective planes (\cite[Theorem
  3.6]{Ertl-Hubl}), any smooth curve in a weighted projective plane must
  be irreducible; indeed any two irreducible components must intersect
  nontrivially, and the intersection lies in the singular locus of the
  curve. Therefore, $X$ is irreducible.
  
  The genus formula can be proved by applying the Riemann--Hurwitz theorem to
  the morphism $Y_{\Q} \to X$ from the curve
  $Y_{\Q}\colon c_1\y_1^{\ell} + c_2\y_2^{\ell} + c_3\y_3^{\ell} = 0 \subset \P^2_{\Q}$ given
  by
  \[(y_1:y_2:y_3) \mapsto (y_1^{\ell/e_1}:y_2^{\ell/e_2}:y_3^{\ell/e_3}),\] for
  $\ell = \lcm(e_1,e_2,e_3)$. See \cite[Corollary A.8]{KobinKeyes26}.
\end{proof}

\subsection{The Belyi cover}
\label{sec:belyi-cover}
In this section, we will study the properties of the $\Q$-morphism
\begin{equation}
  \label{eq:phi-belyi}
  \phi\colon X \to \Pone_{\Q}, \quad (x_1:x_2:x_3) \mapsto (-c_1x_1^{e_1}:c_3x_3^{e_3}).
\end{equation}

\begin{remark}
  For a morphism of $S$-schemes $\phi\colon X \to Y$, we let $\Autsch(\phi)$ denote the
  \cdef{automorphism $S$-group scheme} of $\phi$. That is, for any $T\to S$,
  $\Autsch(\phi)(T) \colonequals \Aut(\phi_T)$ is the subgroup of $g \in \Aut(X_T)$ which
  satisfies $\phi_T \circ g = \phi_T$.
\end{remark}

Our first task will be to understand $\Autsch(\phi)$ for $\phi$ given by
\cref{eq:phi-belyi}. To this end, consider
$\mu_{e_1} \times \mu_{e_2} \times \mu_{e_3}$. It is a finite flat commutative group scheme
over $S = \Spec \Z$ of order $e_1e_2e_3$, acting on $\P(\bfw)$ via coordinatewise
multiplication
\[ (x_1:x_2:x_3)\cdot (\xi_1,\xi_2,\xi_3) = (\xi_1x_1 : \xi_2x_2 : \xi_3x_3).\]

Consider the group homomorphism
\begin{equation}
\label{eqn:iota}
\iota\colon \mu_\ell \to \mathbb{G}_m^{3}, \quad  \xi \mapsto (\xi^{w_1}, \xi^{w_2}, \xi^{w_3}).
\end{equation}

\begin{lemma}
  \label[lemma]{lemma:K-is-the-kernel}
  The map $\iota$
  is a closed immersion with image 
  $\sfK = (\mu_{e_1}\times \mu_{e_2} \times \mu_{e_3}) \cap \Gw$. Moreover, the action of
  $(\mu_{e_1}\times \mu_{e_2} \times \mu_{e_3})_{\Q}$ on $\P(\bfw)_{\Q}$ preserves
  $X$ and commutes with $\phi$, and the kernel of the resulting map
  $(\mu_{e_1}\times \mu_{e_2} \times \mu_{e_3})_{\Q} \to \Autsch(\phi)$ is exactly $\sfK_{\Q}$.
\end{lemma}

\begin{proof}
  The homomorphism $\iota$ lands in $\mu_{e_1}\times \mu_{e_2} \times \mu_{e_3}$ because
  $(\xi^{w_i})^{e_i} = \xi^{\ell} = 1$. It is a monomorphism because
  $\gcd(\bfw) = 1$: choosing $n_i \in \Z$ satisfying $\sum_i n_i w_i = 1$ recovers the identity
  $\xi = \prod_i (\xi^{w_i})^{n_i}$ functorially.

  Since $\xi_i^{e_i} = 1$, the action of
  $\mu_{e_1}\times \mu_{e_2} \times \mu_{e_3}$ preserves the equation
  $\sum_i c_i \x_i^{e_i} = 0$, so it preserves $X$. Since
  $\xi_1^{e_1} = \xi_3^{e_3} = 1$ it fixes $(-c_1\x_1^{e_1} : c_3\x_3^{e_3})$, so
  it commutes with $\phi$.

  A point $(\xi_1, \xi_2, \xi_3)$ of
  $\mu_{e_1}\times \mu_{e_2} \times \mu_{e_3}$ acts trivially on $\P(\bfw)$ if and only if it
  lies in $\Gw$, say $\xi_i = t^{w_i}$. Then
  $t^{\ell} = t^{e_1 w_1} = \xi_1^{e_1} = 1$, so $t \in \mu_\ell$ and
  $(\xi_1, \xi_2, \xi_3) \in \iota(\mu_\ell)$. Conversely
  $\iota(\mu_\ell) \subseteq \Gw$ by construction. Hence the kernel of
  $\mu_{e_1}\times \mu_{e_2} \times \mu_{e_3} \to \Autsch(\phi)$ is
  $(\mu_{e_1}\times \mu_{e_2} \times \mu_{e_3}) \cap \Gw = \sfK$.
\end{proof}

\begin{lemma}
  \label[lemma]{lemma:degree}
  The map $\phi$ has degree $e_1e_2e_3/\ell$, and the quotient group
  \[(\mu_{e_1}\times \mu_{e_2} \times \mu_{e_3})(\Qbar)/\sfK(\Qbar)\] acts simply transitively
  on every geometric fiber of $\phi$ away from $\{0,1,\infty\}$.
\end{lemma}

\begin{proof}
  Work over $\Qbar$ and put $\sfu_i \colonequals c_i\x_i^{e_i}$, each of degree
  $\ell$. Each $(x_1, x_2, x_3) \in X(\Qbar)$ lying over
  $(s:t) \in \Pone(\Qbar) \setminus \{0,1,\infty\}$ satisfies
  \[
    u_1 + u_2 + u_3 = 0 \qquad\text{and}\qquad t\,u_1 + s\,u_3 = 0,
  \]
  which determines the point
  $(u_1 : u_2 : u_3) \in \P^2(\Qbar)$ uniquely, with all three coordinates
  nonzero. Fix a point $P_0$ of the fiber and a representative of it in the
  punctured cone. Given any other point $P$ of the fiber, scale a cone
  representative of $P$ by the $t \in \Gw(\Qbar)$ with $t^{\ell}$ equal to the
  relevant ratio, so that $u_i(P) = u_i(P_0)$ for every position $i$. Then
  $\x_i(P)^{e_i} = \x_i(P_0)^{e_i}$, so $\x_i(P) = \xi_i \x_i(P_0)$ for some
  $\xi_i \in \mu_{e_i}$, i.e.\ $P =  P_0\cdot\vec \xi$ with
  $\vec \xi \in (\mu_{e_1}\times \mu_{e_2} \times \mu_{e_3})(\Qbar)$. So
  $(\mu_{e_1}\times \mu_{e_2} \times \mu_{e_3})(\Qbar)$ acts transitively on the fiber. The
  stabilizer of $P_0$ is $\sfK(\Qbar)$ by \Cref{lemma:K-is-the-kernel}, all
  coordinates of $P_0$ being nonzero. Therefore the general fiber has
  $\#(\mu_{e_1}\times \mu_{e_2} \times \mu_{e_3})(\Qbar)/\#\sfK(\Qbar) = e_1e_2e_3/\ell$ points.
\end{proof}

\begin{proposition}
  \label[proposition]{prop:aut}
  The morphism $\phi_\Q \colon X \to \Pone_\Q$ is a geometrically Galois Belyi
  map, and the scaling action induces an isomorphism of $\Q$-group schemes
  \begin{equation}
    \label{eq:aut}
    (\mu_{e_1}\times \mu_{e_2} \times \mu_{e_3})_{\Q}/\sfK_{\Q} \;\xrightarrow{\ \sim\ }\; \Autsch(\phi_\Q), \qquad
    \#\Aut(\phi_{\Qbar}) = \deg \phi_\Q = e_1e_2e_3/\ell.
  \end{equation}
\end{proposition}

\begin{proof}
  Let $G \colonequals (\mu_{e_1}\times \mu_{e_2} \times \mu_{e_3})$. By
  \Cref{lemma:K-is-the-kernel} the scaling action gives an injection
  $G/\sfK \hookrightarrow \Autsch(\phi)$, so
  $\#\Aut(\phi_{\Qbar}) \geq e_1e_2e_3/\ell = \deg\phi_\Q$ by \Cref{lemma:degree}. On the
  other hand a nontrivial automorphism of a curve has finitely many fixed
  points, so choosing a fiber of $\phi_{\Qbar}$ of cardinality $\deg\phi$ and
  disjoint from all of them, the group $\Aut(\phi_{\Qbar})$ acts freely on it and
  $\#\Aut(\phi_{\Qbar}) \leq \deg\phi_\Q$. Hence equality holds throughout: the cover is
  Galois and \cref{eq:aut} is an isomorphism. Everything is defined over $\Q$,
  so the isomorphism is one of $\Q$-group schemes.

  For the Belyi property: $\phi$ is unramified away from $\x_1\x_2\x_3 = 0$ by
  \Cref{lemma:degree}, and the three loci $\x_1 = 0$, $\x_2 = 0$, $\x_3 = 0$
  map respectively to $(0:1)$, $(1:1)$ and $(1:0)$ (the middle one because
  $\x_2 = 0$ forces $-c_1\x_1^{e_1} = c_3\x_3^{e_3}$). 
\end{proof}

For the counting problem it is the \emph{rational} automorphisms that matter.

\begin{lemma}
  \label{lemma:aut-rational}
  Suppose that $e_{i} = \ell$ for some index $i$, or equivalently that
  $w_{i} = 1$. Then
  \[
    \Autsch(\phi) \;\cong\; \prod_{j \neq i} \mu_{e_j}
    \qquad\text{and}\qquad
    \#\Aut(\phi_{\Q}) = \prod_{j \neq i} \gcd(e_j, 2).
  \]
  Explicitly, every $\Q$-rational automorphism of $\phi$ is a sign change
  $\x_j \mapsto \pm\x_j$ for $j \neq i$.
\end{lemma}

\begin{proof}
  Assume, without loss of generality, that $e_3 = \ell$ (or equivalently
  $w_3=1$). Write $G \colonequals \mu_{e_1}\times \mu_{e_2} \times \mu_{e_3}$ and define
  \[
    \pi\colon G \to \mu_{e_1}\times \mu_{e_2}, \qquad
    (\xi_1, \xi_2, \xi_3) \mapsto \paren{\xi_1 \xi_{3}^{-w_1}, \xi_2 \xi_{3}^{-w_2}}.
  \]
  This lands in $\mu_{e_1}\times \mu_{e_2}$, since
  $\paren{\xi_j\xi_{3}^{-w_j}}^{e_j} = \xi_{3}^{-\ell} = 1$ using that
  $e_3 = \ell$. It is a homomorphism of group schemes over $\Z$, and it is split
  surjective by the section sending $(\eta_j)_j$ to the tuple with third
  coordinate $1$. Its kernel consists of the tuples with
  $\xi_j = \xi_{3}^{w_j}$ for every $j$ (including $j = 3$, as $w_3 = 1$). By
  \Cref{lemma:K-is-the-kernel}, that means that $\ker \pi = \sfK$. Now apply
  \Cref{prop:aut}. Finally $\#\mu_n(\Q) = \gcd(n,2)$.
\end{proof}

\subsection{Euclidean signatures}
\label{sec:euclidean-signatures}
We say that $\bfe$ is \cdef{Euclidean} if $\chi = 0$. There are only three
possibilities for the multiset $\{e_1,e_2,e_3\}$.
\tessfigure

\begin{lemma}
  \label[lemma]{lemma:euclidean-signatures}
  A signature $\bfe$ is Euclidean if and only if the multiset $\{e_1,e_2,e_3\}$
  is one of $\{3,3,3\}, \{3,2,6\},$ or $\{4,2,4\}$.
\end{lemma}
  Indeed, these signatures are in bijective correspondence with the triangular
  tessellations of the Euclidean plane. See \cite[Section II.4]{Magnus74} for
  classical background on Euclidean triangle groups, and
  \cite{Radosevich&Voight-23} for a modern exposition of the subject with an
  emphasis on computing the associated Belyi maps.

\begin{remark}
  \label{remark:b=2}
  When $\bfe$ is Euclidean, we will use the convention that
  $e_2 \leq e_1 \leq e_3$. If $\bfe \neq (3,3,3)$, then $e_2 = 2$.
\end{remark}
The weights and bad locus of Euclidean signatures, in accordance with
\cref{eq:bad} are given in \Cref{table:euclidean-weights}. We see that every Euclidean signature satisfies the hypothesis $\ell = e_3$ of~\Cref{lemma:aut-rational}. 

\begin{table}[htbp]
  \centering
  \caption{Euclidean data for $\ggfe = 0$.}
  \label{table:euclidean-weights}
  \rowcolors{2}{white}{gray!15}
  \begin{tabular}{cccc}
    \toprule
    $\bfe$ & $\bfw$ & $w_i > 1$ & $\bad_F$ \\
    \midrule
    $(3,3,3)$ & $(1,1,1)$ & --- & $1$ \\
    $(3,2,6)$ & $(2,3,1)$ & $w_1, w_2$ & $c_1c_2$ \\
    $(4,2,4)$ & $(1,2,1)$ & $w_2$ & $c_2$ \\
    \bottomrule
  \end{tabular}
\end{table}

\begin{table}[htbp]
  \centering
  \caption{Rational automorphisms of $\phi\colon (x_1:x_2:x_3) \mapsto
    (-c_1x_1^{e_1}:c_3x_3^{e_3})$ for the Euclidean signatures.}
  \label{table:euclidean-rational-auts}
  \rowcolors{2}{white}{gray!15}
  \begin{tabular}{ccccl}
    \toprule
    $\bfe$ & $\Autsch(\phi)$ & $\deg\phi $ & $\#\Aut(\phi_{\Q})$ & Rational automorphisms \\
    \midrule
    $(3,3,3)$ & $\mu_3\times\mu_3$ & $9$ & $1$ & The identity \\
    $(3,2,6)$ & $\mu_3\times\mu_2$ & $6$ & $2$ & $(\x_1:\x_2:\x_3)\mapsto(\x_1:\pm\x_2:\x_3)$ \\
    $(4,2,4)$ & $\mu_2\times\mu_4$ & $8$ & $4$ & $(\x_1:\x_2:\x_3)\mapsto(\x_1:\pm\x_2:\pm\x_3)$ \\
    \bottomrule
  \end{tabular}
\end{table}

% \begin{remark}
%   \label{remark:aut-rational-independent}
%   The order of $\Aut(\phi_{\Q})$ depends only on the signature $\bfe$, and not on
%   the coefficients $c_1,c_2,c_3$: the automorphisms are scalings, which do not
%   involve the coefficients. Note also that one must compute in $(G/\sfK)(\Q)$
%   rather than in $G(\Q)/\sfK(\Q)$, which can be too small: in signature
%   $\{4,2,4\}$ the class of the irrational tuple $(i,1,i)$ equals the class of
%   $(-1,-1,-1)$, both normalizing to $(1,-1,1)$.
% \end{remark}

\section{Fermat curves of genus one}
\label{sec:fermat-genus-1} In this section, we write down plane models of the genus one
Fermat curves and their Jacobians. We keep the coordinates $\x_1,\x_2,\x_3$ on
$\P(\bfw)$ from \Cref{sec:generalities}, and reserve $\x,\y,\z$ for coordinates on
$\P^2$.

\subsection{Plane models}
\label{sec:plane-models}
Suppose that $\bfe$ is Euclidean, with the convention $e_2 \leq e_1 \leq e_3$ of
\Cref{remark:b=2}. A smooth plane curve of genus one is a cubic, so $X$ has a smooth
plane model over $\Q$ if and only if it has a very ample $\Q$-rational divisor of
degree $3$.

If $\bfe = (3,3,3)$, then $\bfw = (1,1,1)$, and $X \subset \P(\bfw) = \P^2$ is
already the smooth plane cubic $c_1\x_1^3 + c_2\x_2^3 + c_3\x_3^3 = 0$.

If $\bfe = (3,2,6)$, then $\bfw = (2,3,1)$. Let $C \subset \P^2_\Q$ be the cubic
$c_1\x^3 + c_2\y^2\z + c_3\z^3 = 0$. Then
\[
  X \to C, \qquad (\x_1:\x_2:\x_3) \mapsto (\x_1\x_3 : \x_2 : \x_3^3)
\]
is an isomorphism. Indeed, its coordinates all have degree $3$; they vanish
simultaneously on $X$ only if $\x_2 = \x_3 = 0$, and hence $\x_1 = 0$; and the
equation of $C$ pulls back to $\x_3^3$ times the equation of $X$. So it is a
morphism. It restricts to an isomorphism from the chart $\x_3 \neq 0$ onto the chart
$\z \neq 0$, with inverse $(\x:\y:\z) \mapsto (\x\z : \y\z^2 : \z)$. The only point
of $C$ with $\z = 0$ is $(0:1:0)$, where the partial derivative of the equation of $C$
with respect to $\z$ is $c_2 \neq 0$; hence $C$ is smooth, and a birational morphism
of smooth projective curves is an isomorphism.

If $\bfe = (4,2,4)$, then $\bfw = (1,2,1)$. The weighted projective curve
$X$ is smooth, but it might not admit a smooth plane model: $X$ has a
smooth plane model over $\Q$ if and only if $X(\Q) \neq \emptyset$. Indeed, the section
$\x_3 = 0$ is a $\Q$-rational divisor of degree $2$ on $X$. Subtracting it
from a $\Q$-rational divisor of degree $3$ would give one of degree $1$; by
Riemann--Roch, its class contains a unique effective divisor, which is then a
rational point. Conversely, if $P_0 \in X(\Q)$, then $|3P_0|$ embeds $X$
as a plane cubic. Since $X(\Q)$ may be empty, we work with
$X \subset \P(1,2,1)$ directly in this case.

\subsection{\texorpdfstring{$j$}{j}-invariants} Suppose that $\bfe$ is Euclidean, and
observe that under the convention $e_2 \leq e_1 \leq e_3$ of \Cref{remark:b=2}, we
have $e_3 = 3, 6, 4$ across the three signatures in
\Cref{table:euclidean-weights}.

\begin{lemma}
  \label[lemma]{lemma:j-invariant} Suppose that $\bfe$ is Euclidean. Let
  $\zeta$ be a primitive $e_3$-th root of unity. Then
  \[
    [\zeta] \colon X_{\Qbar} \to X_{\Qbar}, \qquad
    (x_1:x_2:x_3) \mapsto ( x_1 : x_2 : \zeta x_3)
  \]
  is an automorphism of order $e_3$ with a fixed point. Consequently
  $\Aut(\Jac(X_{\Qbar}), O)$ contains an element of order $e_3$, and
  \[
    j(\Jac(X)) =
    \begin{cases}
      0 & \text{if } (e_1,e_2,e_3) = (3,3,3), (3,2,6), \\
      1728 & \text{if } (e_1,e_2,e_3) = (4,2,4).
    \end{cases}
  \]
\end{lemma}

\begin{proof}
  The map $[\zeta]$ preserves the equation because $\zeta^{e_3} = 1$. Moreover,
  $[\zeta]$ is nontrivial on $\P(\bfw)$: otherwise, it would be induced by some
  $t \in \Gw$ satisfying $t^{w_3} = \zeta$ and $t^{w_1} = t^{w_2} = 1$, hence
  $t = 1$ by $\gcd(w_1,w_2) = 1$, contradicting $\zeta \neq 1$. The same computation
  shows $[\zeta]$ has order $e_3$. If $(x_1,x_2,x_3) \in X(\Qbar)$ is fixed by
  $[\zeta]$, then $x_3 = 0$. Since the equation
  $c_1\x_1^{e_1} + c_2\x_2^{e_2} = 0$ has $\Qbar$-points, we can take one such
  point $O$ as the origin of $X_{\Qbar} \cong \Jac(X_{\Qbar})$. This choice
  realizes $\langle[\zeta]\rangle$ in $\Aut(\Jac(X_{\Qbar}), O)$ as a cyclic subgroup of order
  $3$, $6$ or $4$. It has an element of order $3$ only when $j = 0$ and one of
  order $4$ only when $j = 1728$. Finally, since $j$ is rational,
  $j(\Jac(X)) = j(\Jac(X_\Qbar))$.
\end{proof}

\subsection{Signature \texorpdfstring{(3,3,3)}{3,3,3}}

Here $\bfw = (1,1,1)$, so $\P(\bfw) = \P^2$ and $X$ is a plane cubic. It need not have a rational point,
even if it has $\Q_\nu$-points for every place $\nu$, as witnessed by Selmer's
famous example
\[ 3\x_1^3 + 4\x_2^3 + 5\x_3^3 = 0.\]

Even when $X(\Q) = \emptyset$, the Jacobian $\Jac(X)$ is an elliptic curve over
$\Q$; the following gives explicit models for it.

\begin{proposition}
  \label[proposition]{prop:333-jacobian}
  Suppose that $\bfe = (3,3,3)$. Then $\Jac(X)$ is isomorphic over $\Q$ to the
  plane cubic
  \[
    J \colon \x^3 + \y^3 = c_1c_2c_3\,\z^3
  \]
  with origin $(1:-1:0)$, and to the elliptic curve in Weierstrass form
  \[
    E \colon \y^2\z = \x^3 - 432(c_1c_2c_3)^2\z^3
  \]
  with origin $(0:1:0)$.
\end{proposition}

\begin{proof}
  The plane cubic $J$ is a model of $\Jac(X)$ by \cite[Theorem~1]{Selmer-51}.
  Its rational point $(1:-1:0)$ is a flex, which we choose as the origin. The maps
  \begin{align*}
    J \to E, \qquad (\x:\y:\z) &\mapsto
      \paren{12c_1c_2c_3\,\z : 36c_1c_2c_3\,(\x - \y) : \x + \y}, \\
    E \to J, \qquad (\x:\y:\z) &\mapsto
      \paren{\y + 36c_1c_2c_3\,\z : -\y + 36c_1c_2c_3\,\z : 6\x}
  \end{align*}
  are morphisms that respect the origins $(1:-1:0)$ and $(0:1:0)$, and they are
  mutually inverse (their composition is multiplication of the coordinates by
  $72c_1c_2c_3$). Hence they are $\Q$-isomorphisms of elliptic curves. 
  % Alternatively, one can identify $E$ directly from the invariants of $X$. In
  % characteristic not $2$ or $3$, the Jacobian of a genus one model with
  % \cdef{geometric invariants} $c_4, c_6$ is
  % $\y^2\z = \x^3 - 27c_4\x\z^2 - 54c_6\z^3$
  % \cite[Proposition~2.3]{Fisher06}.
  % \begin{magma}
  %   F<[c]> := FieldOfFractions(PolynomialRing(Rationals(),3));
  %   Pol<[x]> := PolynomialRing(F,3);
  %   eqn := c[1]*x[1]^3 + c[2]*x[2]^3 + c[3]*x[3]^3;
  %   PP2 := ProjectiveSpace(F, 2);
  %   X := Curve(PP2,eqn);
  %   cInvariants(GenusOneModel(X)); // Gives (*$[c_4, c_6] = [0,\, 2^33^6(c_1c_2c_3)^2].$*)
  % \end{magma}
  % For the ternary cubic defining $X$ we have
  % $c_4(X) = 0$ and $c_6(X) = 2^33^6(c_1c_2c_3)^2$, and
  % $(\x:\y:\z) \mapsto (9\x : 27\y : \z)$ is an isomorphism from $E$ to this curve.
  % Finally, $X$ and $J$ become isomorphic over $\Qbar$. Choose
  % $\gamma_i \in \Qbar$ such that $\gamma_i^3 = c_i$. Then
  % \[
  %   \psi\colon X_{\Qbar} \to J_{\Qbar}, \qquad
  %   (\x_1:\x_2:\x_3) \mapsto (\gamma_1^2\gamma_2\,\x_1 : \gamma_1\gamma_2^2\,\x_2 : -\x_3)
  % \]
  % is a linear change of coordinates that pulls the equation of $J$ back to $c_1c_2$
  % times that of $X$. It is therefore an isomorphism, and it sends the flex
  % $(\gamma_2:-\gamma_1:0)$ of $X_{\Qbar}$ to the origin $(1:-1:0)$ of $J$.
\end{proof}

\subsection{Signature \texorpdfstring{(3,2,6)}{3,2,6}}

In contrast to the other signatures, when $\bfe = (3,2,6)$ the genus one curve
$X$ always has a rational point!

\begin{proposition}
  \label[proposition]{prop:326}
  The genus one curve $X$ has the rational point $P_0 \colonequals (-c_1c_2
  : c_1^2c_2 : 0)$. Choosing it as the identity, the elliptic curve $X$ is
  isomorphic to the short Weierstrass model
  \[E \colon \y^2\z = \x^3 - c_1^2c_2^3c_3\z^3 \]
  via
  \[
    \psi \colon X \longrightarrow E,
    \qquad
    (\x_1:\x_2:\x_3) \mapsto (-c_1c_2\x_1\x_3 : c_1c_2^2\x_2 : \x_3^3).
  \]
\end{proposition}

\begin{proof}
  At $\x_3 = 0$ the defining equation reads $c_1\x_1^3 + c_2\x_2^2 = 0$, and
  \[
    c_1(-c_1c_2)^3 + c_2(c_1^2c_2)^2 = -c_1^4c_2^3 + c_1^4c_2^3 = 0,
  \]
  so $P_0 \in X(\Q)$. Since $\bfw = (2,3,1)$, the three coordinates of $\psi$ have
  degrees $2+1$, $3$ and $1+1+1$, so $\psi$ is a well-defined morphism to $\P^2$.
  Using $c_2\x_2^2 = -(c_1\x_1^3 + c_3\x_3^6)$, we see that
  \begin{align*}
    \x^3 - c_1^2c_2^3c_3\z^3
    &= -c_1^3c_2^3\,\x_1^3\x_3^3 - c_1^2c_2^3c_3\,\x_3^9 \\
    &= -c_1^2c_2^3\,\x_3^3\paren{c_1\x_1^3 + c_3\x_3^6} = c_1^2c_2^4\,\x_2^2\x_3^3 = \y^2\z,
  \end{align*}
  so $\psi$ lands on $E$. It carries $P_0$ to
  $(0 : c_1^3c_2^3 : 0) = (0:1:0)$, and $E$ is smooth because
  $c_1c_2c_3 \neq 0$; because a birational map of smooth projective curves is an isomorphism, $\psi$ is an isomorphism.
\end{proof}

\subsection{Signature \texorpdfstring{(4,2,4)}{4,2,4}}

Here $\bfw = (1,2,1)$. As in $\{3,3,3\}$, $X$ need not have a rational
point.

\begin{proposition}
  \label[proposition]{prop:424-jacobian}
  The genus one curve $X$ has Jacobian isomorphic to the Weierstrass model
  \[E\colon \y^2\z = \x^3 - 4c_1c_2^2c_3\,\x\z^2.\]
\end{proposition}

\begin{proof}
  We compute the geometric invariants $c_4(X)$ and $c_6(X)$ in
  \Magma~\cite{Magma}.
  \begin{magma}
  F<[c]> := FieldOfFractions(PolynomialRing(Rationals(),3));
  Pol<[x]> := PolynomialRing(F,3);
  eqn := c[1]*x[1]^4 + c[2]*x[2]^2 + c[3]*x[3]^4;
  PPw := ProjectiveSpace(F, [1,2,1]);
  X := Curve(PPw,eqn);
  cInvariants(GenusOneModel(X)); // Gives (*$[c_{4}, c_6] = [192\, c_1c_3/c_2^2, 0].$*)
\end{magma}
  Scaling by $\lambda = c_2/2$, we conclude that
  $c_4(X) = 12\, c_1c_2^2c_3$ and $c_6(X) = 0$, and the elliptic curve
  $E$ above is visibly isomorphic to $\y^2\z = \x^3 - 27\,c_4(X)\,\x\z^2 -
  54\,c_6(X)\,\z^3$. The result follows again from \cite[Proposition 2.3]{Fisher06}.
\end{proof}

% A geometric isomorphism $\psi\colon X_{\Qbar}\to E_{\Qbar}$ can be obtained as
% follows. Choose $\gamma_1, \gamma_2, \gamma_3 \in \Qbar$ with $\gamma_1^4 = c_1$,
% $\gamma_2^2 = -c_2$ and $\gamma_3^4 = c_3$, and put
% $\delta \colonequals \gamma_1\gamma_2\gamma_3$, so that $\delta^4 = c_1c_2^2c_3$.
% Take $P_0 \colonequals (\gamma_2:\gamma_1^2\gamma_2:0)$ as the origin of $X_{\Qbar}$.
% Consider the curves $Y \colon \y_2^2 = \y_1^4 + \y_3^4 \subset \P(1,2,1)_{\Qbar}$ and
% $E_0 \colon \y^2\z = \x^3 - 4\x\z^2$. The maps
% \begin{align*}
%   X_{\Qbar} \to Y \colon (x_1:x_2:x_3) &\mapsto (\gamma_1 x_1 : \gamma_2 x_2 : \gamma_3 x_3), \\
%   Y \to E_0 \colon (y_1:y_2:y_3) &\mapsto (2y_3^2 : 4y_1y_3 : y_2 - y_1^2), \\
%   E_0 \to E_{\Qbar} \colon (x:y:z) &\mapsto (\delta^2 x : \delta^3 y : z)
% \end{align*}
% are isomorphisms: the middle map is defined away from $(1:1:0)$, and since
% $(y_2 - y_1^2)(y_2 + y_1^2) = y_3^4$ on $Y$, it agrees with
% $(2(y_1^2+y_2)y_3 : 4y_1(y_1^2+y_2) : y_3^3)$, which sends $(1:1:0)$ to $(0:1:0)$.
% Their composition is the $\Qbar$-isomorphism of elliptic curves
% \begin{equation*}
%   \psi\colon X_{\Qbar} \to E_{\Qbar}, \quad (x_1:x_2:x_3) \mapsto
%   (-2c_2c_3\gamma_1^2\, x_3^2 : -4c_1c_2c_3\gamma_2\, x_1x_3 : \gamma_2x_2 - \gamma_1^2x_1^2),
% \end{equation*}
% which sends $P_0$ to $(0:1:0)$.

\subsection{Summary} We summarize the results of this section in the following table.
\begin{table}[htbp]
  \centering
  \caption{Summary of \Cref{sec:fermat-genus-1}.}
  \rowcolors{2}{white}{gray!15}
  \begin{tabular}{clcl}
    \toprule
    $\bfe$ & $\Jac(X) \subseteq \P^2$ & $j$-invariant & $X(\Q) \neq \emptyset$\,? \\
    \midrule
    $(3,3,3)$ & $\y^2\z = \x^3 - 432(c_1c_2c_3)^2\z^3$ & $0$ & not in general \\
    $(3,2,6)$ & $\y^2\z = \x^3 - c_1^2c_2^3c_3\,\z^3$ & $0$ & always \\
    $(4,2,4)$ & $\y^2\z = \x^3 - 4c_1c_2^2c_3\,\x\z^2$ & $1728$ & not in general \\
    \bottomrule
  \end{tabular}
  \label{tab:jacobians}
\end{table}

\section{Counting rational points in the image of a map from an elliptic curve}
\label{sec:counting}

Throughout this section, we fix the following notation.
\begin{itemize}[leftmargin=*]
  \item $E/\Q$ is an elliptic curve of Mordell--Weil rank $r$.
  \item $\phi\colon E\to\Pone_\Q$ is a finite map of degree
    $d\colonequals\deg\phi\geq 2$.
  \item $\mcL\colonequals\phi^*\OO_{\Pone}(1)$ is the pullback line bundle,
    which has degree $d$ on $E$.
  \item $\widehat\logHt_\mcL$ is the canonical height attached to $\mcL$.
  \item $\langle-,-\rangle_\mcL$ is the inner product on $E(\Q)\otimes_\Z\R$
    determined by
    \[2\langle x,y\rangle_\mcL\colonequals\widehat\logHt_\mcL(x+y)
      -\widehat\logHt_\mcL(x)-\widehat\logHt_\mcL(y)\] for $x,y\in E(\Q)$.
  \item $\Reg_\mcL(E/\Q)\colonequals\det\paren{\langle x_i,x_j\rangle_\mcL}_{i,j}$
    is the \cdef{regulator of $E/\Q$ with respect to $\mcL$}, where
    $x_1,\dots,x_r$ is a $\Z$-basis of $E(\Q)/E(\Q)_{\tors}$.
  \item $\Reg(E/\Q)\colonequals\Reg_{\OO_E(2O)}(E/\Q)$, where $O$ is the origin
    of $E$. This is the normalization of the regulator used by Magma and the
    LMFDB (\Cref{remark:regulator-normalisation}).
  \item For a subgroup $\Lambda\le E(\Q)$ of finite index and a point
    $P_0\in E(\Q)$, the \cdef{counting function} of $\phi$ on the coset
    $P_0+\Lambda$ is
    \begin{equation}
      \label{def:N_phi(h)}
      N_\phi(H;P_0+\Lambda)\colonequals
      \#\brk{P\in P_0+\Lambda : \Ht(\phi(P))\leq H}.
    \end{equation}
  \item For a subgroup $\Lambda\le E(\Q)$ of finite index, we set
    \begin{equation}
      \label{eqn:delta-phi}
      \delta(\phi;\Lambda)\colonequals\frac{\pi^{r/2}}{\Gamma(1+r/2)}\cdot
      \frac{\#E(\Q)_{\tors}}{[E(\Q):\Lambda]\sqrt{\Reg_\mcL(E/\Q)}}>0,
    \end{equation}
    and $\delta(\phi)\colonequals\delta(\phi;E(\Q))$.
  \item The number of rational points of bounded height in the image of
    $\phi$ is
    \begin{equation}
      \label{def:N(phi(E(Q));h)}
      N(\phi(E(\Q));H)\colonequals
      \#\brk{Q\in\phi(E(\Q)) : \Ht(Q)\leq H}.
    \end{equation}
  \end{itemize}
  
  We will need the following consequence of a classical theorem of N\'eron
  \cite[Theorem B.6.3]{Hindry&Silverman00},
  \cite[p.~54]{Serre97-LecturesMordellWeil}.

\begin{proposition}
  \label{prop:neron-corollary}
  Let $\Lambda\le E(\Q)$ be a subgroup of finite index, and let $P_0\in E(\Q)$.
  Then
  \[
    N_\phi(H;P_0+\Lambda)=\delta(\phi;\Lambda)\cdot(\log H)^{r/2}
    +O\!\paren{(\log H)^{(r-1)/2}},\quad\text{as }H\to\infty,
  \]
  where the implied constant depends on $\phi$ and $\Lambda$.
\end{proposition}

\begin{proof}
    We work throughout this proof with the logarithmic height
    $\logHt_{\Pone(\Q)}\colonequals \log\circ\Ht$ as opposed to the multiplicative height
    appearing in the statement. With respect to this height function, the
    statement is equivalently that
    \[\#\brk{P\in P_0+\Lambda : \logHt_{\Pone(\Q)}(\phi(P))\le h}=\delta(\phi;\Lambda)\cdot
    h^{r/2}+O\!\paren{h^{(r-1)/2}}, \quad\text{as }h\to\infty,\] 
    with $\delta(\phi;\Lambda)$ as in \cref{eqn:delta-phi}.
    Let $\tau\colon E\to E$, $Q\mapsto P_0+Q$, denote translation by $P_0$, and
    let $\mcL'\coloneqq\tau^*\mcL$. Note that $\deg\mcL'=\deg\mcL=d$. It follows from \cite[Theorem B.3.1]{Hindry&Silverman00} that
    $P\mapsto\logHt_{\Pone(\Q)}(\phi(P_0+P))$ is a logarithmic Weil height for
    $\mcL'$. In particular, letting $\widehat\logHt_{\mcL'}$ denote the canonical
    height associated to $\mcL'$, we have that
    \[\left|\widehat\logHt_{\mcL'}-\logHt_{\Pone(\Q)}\circ\phi\circ\tau\right|\le O(1).\] 
    Hence, it suffices to prove the corresponding asymptotic for
    \[N_{\mcL'}(h;\Lambda)\coloneqq\#\brk{P\in\Lambda : \widehat\logHt_{\mcL'}(P)\le h}.\]
    Since $d=\deg\mcL'>0$, N\'eron's theorem \cite[Theorem B.6.3]{Hindry&Silverman00}
    (stated there in terms of a multiplicative height function) shows that there
    is a constant $c>0$ such that
    \[N_{\mcL'}(h;\Lambda)= c\cdot h^{r/2}+O\paren{h^{(r-1)/2}}\] as
    $h\to\infty$. It remains to show that $c=\delta(\phi;\Lambda)$. Let
    $V\coloneqq E(\Q)\otimes_\Z\R\cong\R^r$.
    Let $\langle-,-\rangle_{\mcL'}$ be the inner product on $V$ attached to
    $\widehat\logHt_{\mcL'}$, defined as for $\mcL$ at the start of this section.
    Let $B\subset \R^r$ denote the unit ball, and let
    $L \coloneqq \Lambda/\Lambda_{\tors}\subset V$. Then,
    the proofs of \cite[Theorem B.6.3 and Lemma B.6.4]{Hindry&Silverman00} (see also \cite[p. 54]{Serre97-LecturesMordellWeil}) show that
    \begin{equation}\label{eqn:delta-phi-first-form}
      c=\lambda(B)\cdot\frac{\#\Lambda_{\tors}}{\covol
        L}=\frac{\pi^{r/2}}{\Gamma(1+r/2)}\cdot \dfrac{\#\Lambda_{\tors}}{\covol L},
    \end{equation}
    where $\lambda$ denotes the usual Lebesgue measure on $\R^r$.
    Fix a $\Z$-basis $x_1,\dots,x_r$ of $L$. Then, the covolume of $L$ above is equal to the square root of the determinant of the matrix $\left(\langle{x_i,x_j}\rangle_{\mcL'}\right)_{i,j}$. When $L=L_0\coloneqq E(\Q)/E(\Q)_{\tors}$, this is, by definition, the square root of the regulator $\Reg_{\mcL'}(E/\Q)$. In general, one has $\covol L=[L_0:L]\covol L_0=[L_0:L]\sqrt{\Reg_{\mcL'}(E/\Q)}$. Since $[E(\Q):\Lambda]=[L_0:L][E(\Q)_{\tors}:\Lambda_{\tors}]$, we conclude from \cref{eqn:delta-phi-first-form} that
    \[c=\frac{\pi^{r/2}}{\Gamma(1+r/2)}\cdot\frac{\#E(\Q)_{\tors}}{[E(\Q):\Lambda]\sqrt{\Reg_{\mcL'}(E/\Q)}}.\]
    Comparing with \cref{eqn:delta-phi}, it remains to show that
    $\Reg_{\mcL'}(E/\Q)=\Reg_\mcL(E/\Q)$. Since canonical heights are additive in
    the line bundle, the difference $\langle-,-\rangle_{\mcL'}-\langle-,-\rangle_\mcL$
    is the pairing attached in the same way to
    $\mathcal M\coloneqq\mcL'\otimes\mcL^{-1}$. Since $\deg\mathcal M=0$, the
    canonical height $\widehat\logHt_{\mathcal M}$ is linear, so this pairing is
    zero. Thus $\langle-,-\rangle_{\mcL'}=\langle-,-\rangle_\mcL$, hence
    $\Reg_{\mcL'}(E/\Q)=\Reg_\mcL(E/\Q)$ and $c=\delta(\phi;\Lambda)$, completing
    the proof.
\end{proof}

  \begin{lemma}
    \label[lemma]{lemma:delta-explicit}
    We have
    \begin{equation}
      \label{eq:delta-explicit}
      \delta(\phi) = \#E(\Q)_{\tors} \cdot
      \frac{\paren{2\pi/d}^{r/2}}{\Gamma\paren{1 + \tfrac r2}\cdot\sqrt{\Reg(E/\Q)}}.
    \end{equation}
  \end{lemma}

  \begin{proof}
    On an elliptic curve, the canonical height attached to a symmetric line
    bundle of degree $k$ is $\tfrac k2$ times the one attached to $\OO_E(2O)$.
    The inner product $\langle-,-\rangle_\mcL$ is $\frac12$ times the inner
    product attached to the symmetric line bundle $\mcL\otimes[-1]^*\mcL$ of
    degree $2d$, so it is $\frac12\cdot\frac{2d}2=\frac d2$ times the one
    attached to $\OO_E(2O)$.
    Thus, $\Reg_\mcL(E/\Q)=(d/2)^{r}\Reg(E/\Q)$, from which the claim follows. 
  \end{proof}

\begin{lemma}\label{lem:counting-image-psi}
    Let $\psi\colon C\to\Pone_\Q$ be a geometrically Galois cover of curves, defined over $\Q$.
    Let $S\subset C(\Q)$ be an $\Aut(\psi)$-stable subset of rational points. Then,
    \[N(\psi(S);H)=\frac1{\#\Aut(\psi)}\paren{\#\{P \in S : \Ht(\psi(P)) \leq H\}+O(1)},\quad\text{as }H\to\infty,\]
    where the implied constant depends on the map $\psi$. 
\end{lemma}
\begin{proof}
    If $S$ is finite, the statement is trivially true, so suppose $\#S=\infty$. Let $P,P'\in S$ be points in the same unramified $\psi$-fiber. Since $\psi_{\Qbar}$ is Galois, there is some $\alpha\in\Aut(\psi_{\Qbar})$ with $P'=\alpha(P)$. Moreover, since $P$ is not a ramification point of $\psi$, then $\alpha$ is unique and, for every $\sigma\in\Gal(\Qbar/\Q)$, we have
    \[\sigma(\alpha)(P)=\sigma(\alpha(P))=\sigma(P')=P'=\alpha(P),\]
    so $\sigma(\alpha)=\alpha$ by uniqueness. Since only finitely many points ramify, we deduce that $\alpha\in\Aut(\psi)$ is defined over $\Q$. Hence, the fibers of $S\to\psi(S)$ are $\Aut(\psi)$-orbits, with $\Aut(\psi)$ acting freely away from the ramification locus.
\end{proof}

\begin{remark}
\label{remark:rational-vs-geometric-aut}
    When applying \Cref{lem:counting-image-psi} to $\phi\colon E\to\P^1$, it is the cardinality of $\Aut(\phi_{\Q}) = \Autsch(\phi)(\Q)$ that appears in
    the leading constant of $N(\phi(E(\Q));H)$, not the cardinality of
    $\Aut(\phi_{\Qbar}) = \Autsch(\phi)(\Qbar)$ (which equals $\deg \phi$).
\end{remark}
  
\section{The Primitive Locus}
\label{sec:minimal-solutions}

In the homogeneous case $\bfe = (3,3,3)$ every rational point of $X$ lifts to a
primitive integral solution, and \Cref{prop:neron-corollary} counts $\Omega_F$
directly. For the other Euclidean signatures this fails
(\Cref{ex:sameera's-example-1}). In this section we show that the rational
points lifting to primitive solutions are cut out by finitely many $p$-adic
conditions, each open and closed, and deduce that they form a finite union of
orbits of a finite index subgroup of $\Jac(X)(\Q)$ (\Cref{thm:primitive-subgroup}).

Throughout this section we assume that $\bfw$ is well-formed. This holds for
every Euclidean signature (\Cref{table:euclidean-weights}), and it makes
\Cref{lemma:U->X-Gm-torsor} available.

\subsection{Minimal lifts}

\begin{definition}
  \label{def:w-minimal}
  Let $p$ be a prime. A triple $(x_1,x_2,x_3) \in \Z_p^3$ is \cdef{$\bfw$-minimal
    at $p$} if $v_p(x_i) < w_i$ for some $i$, that is, if $p^{w_i} \mid x_i$ fails
  for some $i$. A triple in $\Z^3$ is \cdef{$\bfw$-minimal} if it is
  $\bfw$-minimal at every prime. For $A \in \{\Z, \Z_p\}$ we write
  \begin{equation}
    \label{eq:F-Zp}
    \mcF(A) \colonequals \brk{(x_1,x_2,x_3) \in A^3 \;:\; F(x_1,x_2,x_3) = 0, \;
      (x_1,x_2,x_3) \text{ is $\bfw$-minimal}}
  \end{equation}
  for the set of $\bfw$-minimal solutions.
\end{definition}

A $\bfw$-minimal triple is nonzero. Recall that $\mcU(A)$ consists of the
solutions in $A^3$ whose coordinates generate the unit ideal: the primitive
integral solutions when $A = \Z$, and the solutions with
$\min_i v_p(x_i) = 0$ when $A = \Z_p$. Such solutions are $\bfw$-minimal, so
\[
  \mcU(A) \subseteq \mcF(A) \subseteq \mcU(K),
\]
where $K \in \{\Q, \Q_p\}$ is the fraction field of $A$.
The first inclusion is strict in general; see \Cref{ex:sameera's-example-1}.

\begin{lemma}
  \label[lemma]{lemma:minimal-lifts}
  Let $A$ be $\Z$ or $\Z_p$ for a prime $p$, and let $K$ be its field of
  fractions. Then every $P \in X(K)$ admits a lift to $\mcF(A)$, and any two
  lifts differ by the $\bfw$-action of $A^\times$.
\end{lemma}

\begin{proof}
  Since $K$ is an $R$-algebra, \Cref{lemma:U->X-Gm-torsor} shows that
  $\mcU_K \to X_K$ is a $\Gw[,K]$-torsor. As $\HH^1(K,\Gm) = 0$ by Hilbert's
  Theorem~90, $P$ lifts to some $(x_1,x_2,x_3) \in \mcU(K)$. For each prime $q$
  of $A$ (note that $q = p$ if $A = \Z_p$), set $k_q \colonequals -\min_i \lfloor v_q(x_i)/w_i \rfloor$, the
  minimum taken over the $i$ with $x_i \neq 0$. Then $k_q = 0$ for all but
  finitely many $q$. Put $\lambda \colonequals \prod_q q^{k_q} \in K^\times$ and
  $y_i \colonequals \lambda^{w_i}x_i$. For every $q$ and every $i$ with
  $x_i \neq 0$,
  \[
    \lfloor v_q(y_i)/w_i \rfloor = \lfloor v_q(x_i)/w_i \rfloor + k_q \geq 0,
  \]
  with equality for some $i$. Hence, for every $q$, we have $v_q(y_i) \geq 0$
  for all $i$ (trivially if $y_i = x_i = 0$) and $v_q(y_i) < w_i$ for some $i$.
  So $y_i \in A$ for all $i$, and $(y_1,y_2,y_3) \in \mcF(A)$ lies over the point $P$.

  Any two lifts of $P$ lie in $\mcU(K)$ over $P$, so by the torsor property they
  have the form $(x_1,x_2,x_3)$ and
  $(\lambda^{w_1}x_1,\lambda^{w_2}x_2,\lambda^{w_3}x_3)$ for some
  $\lambda \in K^\times$. If both lie in $\mcF(A)$, then for each prime $q$ of
  $A$, $v_q(\lambda) > 0$ contradicts the $\bfw$-minimality of the second triple
  at $q$, and $v_q(\lambda) < 0$ that of the first. Hence $v_q(\lambda) = 0$ for
  every $q$, that is, $\lambda \in A^\times$.
\end{proof}

\begin{definition}
  \label{def:primitive-locus}
  The \cdef{primitive locus} is the image of the projection
  $\pr\colon \mcU(\Z) \to X(\Q)$, $(x_1,x_2,x_3) \mapsto (x_1:x_2:x_3)$:
  \begin{equation}
    \label{eq:primitive-locus}
    X(\Q)_1 \colonequals \pr(\mcU(\Z)).
  \end{equation}
  It consists of the rational points that lift to primitive integral solutions
  of \cref{eq:gfe}. Likewise, for a prime $p$ the set of \cdef{primitive $p$-adic
    points} is
  \begin{equation}
    \label{eq:local-primitive-locus}
    X(\Q_p)_1 \colonequals \pr\paren{\mcU(\Z_p)} \subseteq X(\Q_p).
  \end{equation}
  The \cdef{primitivity defect} of $P \in X(\Q)$ is
  \begin{equation}
    \label{eq:defect}
    \defect(P) \colonequals \gcd(x_1,x_2,x_3),
  \end{equation}
  where $(x_1,x_2,x_3) \in \mcF(\Z)$ is any lift of $P$; by
  \Cref{lemma:minimal-lifts} it does not depend on the lift. Setting
  $X(\Q)_d \colonequals \brk{P \in X(\Q) : \defect(P) = d}$ we obtain a partition
  \begin{equation}
    \label{eq:defect-partition-X}
    X(\Q) = \bigsqcup_{d \geq 1} X(\Q)_d.
  \end{equation}
\end{definition}

The two meanings of $X(\Q)_1$ agree: primitive solutions are $\bfw$-minimal, so
by \Cref{lemma:minimal-lifts} a point $P$ lifts to $\mcU(\Z)$ if and only if its
$\bfw$-minimal lift (unique up to sign) is primitive, that is, if and only if
$\defect(P) = 1$.

\begin{example}
  \label{ex:sameera's-example-1}
  Let $F = \x_1^4 - 8\x_2^2 + \x_3^4$, of signature $(4,2,4)$, so that
  $\bfw = (1,2,1)$ and $\bad_F = -8$. The triple $(2,2,2)$ is a solution, and
  it is $\bfw$-minimal because $2^2 \nmid 2$. By \Cref{lemma:minimal-lifts} the
  only other $\bfw$-minimal lift of $P = (2:2:2)$ is $(-2,2,-2)$, so
  $\defect(P) = 2$ and no primitive integral solution lies over $P$. That the
  prime $2$ divides $c_2 = -8$ is not a coincidence
  (\Cref{prop:global-criterion}).
\end{example}

\subsection{Local conditions}

We show that the primitive locus $X(\Q)_1$ is cut out by the local primitive
loci $X(\Q_p)_1$ at the primes $p \mid \bad_F$.

\begin{lemma}
  \label[lemma]{lemma:local-criterion}
  Let $p$ be a prime.
  \begin{enumerate}[leftmargin=*, label=(\alph*)]
  \item \label{it:local-min-zero} A point $P \in X(\Q_p)$ lies in $X(\Q_p)_1$ if
    and only if $\min_i v_p(x_i) = 0$ for one, equivalently every, lift
    $(x_1,x_2,x_3) \in \mcF(\Z_p)$ of $P$.
  \item \label{it:local-good} If $p \nmid \bad_F$, then $X(\Q_p)_1 = X(\Q_p)$.
  \item \label{it:local-clopen} $X(\Q_p)_1$ is open and closed in $X(\Q_p)$.
  \end{enumerate}
\end{lemma}

\begin{proof}
  \ref{it:local-min-zero} If $\min_i v_p(x_i) = 0$, the lift lies in
  $\mcU(\Z_p)$. Conversely, if $P = \pr(y)$ with $y \in \mcU(\Z_p) \subseteq
  \mcF(\Z_p)$, then every lift in $\mcF(\Z_p)$ is $y$ scaled by a unit
  (\Cref{lemma:minimal-lifts}), which does not affect the valuation.

  \ref{it:local-good} If $p \nmid \bad_F$ then $\Z_p$ is an $R$-algebra, so
  $\mcU_{\Z_p} \to \mcX_{\Z_p}$ is a $\Gw[,\Z_p]$-torsor by
  \Cref{lemma:U->X-Gm-torsor}. Since $\Pic(\Z_p) = 0$, the map
  $\mcU(\Z_p) \to \mcX(\Z_p)$ is surjective, and $\mcX(\Z_p) = X(\Q_p)$ because
  $\mcX$ is proper.

  \ref{it:local-clopen} The set $\mcF(\Z_p)$ is closed in $\Z_p^3$, hence
  compact: it is the zero locus of $F$ minus the open (and closed) subset
  $p^{w_1}\Z_p \times p^{w_2}\Z_p \times p^{w_3}\Z_p$. It is the disjoint union of the
  compact sets $\mcU(\Z_p) = \mcF(\Z_p) \smallsetminus (p\Z_p)^3$ and
  $\mcF(\Z_p) \cap (p\Z_p)^3$. The continuous map $\pr$ sends $\mcF(\Z_p)$ onto
  $X(\Q_p)$ (\Cref{lemma:minimal-lifts}), and by \ref{it:local-min-zero} it
  sends the two pieces onto the disjoint sets $X(\Q_p)_1$ and
  $X(\Q_p) \smallsetminus X(\Q_p)_1$. Both images are compact, hence closed in the Hausdorff
  space $X(\Q_p)$. Since they partition $X(\Q_p)$, each is also open.
\end{proof}

\begin{proposition}
  \label{prop:global-criterion}
  Let $P \in X(\Q)$, let $(x_1,x_2,x_3) \in \mcF(\Z)$ be a lift of $P$, and let
  $p$ be a prime.
  \begin{enumerate}[leftmargin=*, label=(\alph*)]
  \item \label{it:vp-defect} $v_p(\defect(P)) = \min_i v_p(x_i) < \max_i w_i$.
  \item \label{it:p-divides-defect} $p \mid \defect(P)$ if and only if
    $P \notin X(\Q_p)_1$. In particular, $\defect(P)$ divides
    $\prod_{p \mid \bad_F} p^{\max_i w_i - 1}$, so only finitely many defects
    occur.
  \item \label{it:X1-local} The primitive locus is cut out by finitely many local
    conditions:
    \begin{equation}
      \label{eq:X1-local}
      X(\Q)_1 = X(\Q) \cap \bigcap_{p \mid \bad_F} X(\Q_p)_1.
    \end{equation}
  \end{enumerate}
\end{proposition}

\begin{proof}
  \ref{it:vp-defect} The equality is the definition of the $\gcd$. If
  $\min_i v_p(x_i) \geq \max_i w_i$, then $p^{w_i} \mid x_i$ for every $i$,
  contradicting $\bfw$-minimality at $p$.

  \ref{it:p-divides-defect} Since $(x_1,x_2,x_3)$ is $\bfw$-minimal at $p$, it is
  a lift of $P$ in $\mcF(\Z_p)$, and
  \Cref{lemma:local-criterion}\ref{it:local-min-zero} applies. The divisibility
  follows from \ref{it:vp-defect} and
  \Cref{lemma:local-criterion}\ref{it:local-good}.

  \ref{it:X1-local} We saw that $P \in X(\Q)_1$ if and only if $\defect(P) = 1$,
  and by \ref{it:p-divides-defect} this holds if and only if $P \in X(\Q_p)_1$
  for every $p$. By
  \Cref{lemma:local-criterion}\ref{it:local-good}, only the primes
  $p \mid \bad_F$ impose a condition.
\end{proof}

\subsection{The primitive subgroup}

Assume that $\bfe$ is Euclidean, so $X$ is a curve of genus $1$. Let $E\coloneqq\Jac(X)=\Pic^0_{X/\Q}$ be the Jacobian of $X$, an elliptic curve, and note that $E(\Q)$ acts simply transitively on $X(\Q)$. We refer to the action of $E$ on $X$ as \cdef{translation}.

\begin{definition}
  \label{def:primitive-subgroup}
  The \cdef{primitive subgroup} of $E(\Q)$ is the stabilizer of the
  primitive locus under translation,
  \begin{equation}
    \label{eq:primitive-subgroup}
    \Lambda_1 \colonequals \brk{T \in E(\Q) : T + X(\Q)_1 = X(\Q)_1}.
  \end{equation}
\end{definition}

\begin{theorem}[The Primitive Subgroup]
  \label{thm:primitive-subgroup}
  
  Let $\Lambda_1\subset E(\Q)$ be the primitive subgroup. Then:
  \begin{enumerate}[leftmargin=*, label=(\alph*)]
  \item \label{it:Lambda1-finite-index} $\Lambda_1$ has finite index in
    $E(\Q)$.
  \item \label{it:Lambda1-cosets} $X(\Q)_1$ is a finite disjoint union of
    orbits of $\Lambda_1$.
  \end{enumerate}
\end{theorem}
\begin{proof}
    If $X(\Q)_1=\emptyset$, then $\Lambda_1=E(\Q)$ and the statement is clear. Below, assume $X(\Q)_1\neq\emptyset$.
    
    \emph{Part (i): Local subgroups.} Fix $p\mid\bad_F$, and let $\Lambda_p\colonequals\brk{T\in E(\Q_p):T+X(\Q_p)_1=X(\Q_p)_1}$. We claim that $\Lambda_p$ is an open subgroup of $E(\Q_p)$. Note that it suffices to show that $\Lambda_p$ contains some open subgroup of $E(\Q_p)$. For this, we recall, from \Cref{lemma:local-criterion}\ref{it:local-clopen} that $X(\Q_p)_1\subset X(\Q_p)$ is open and closed (hence compact). Consider the map
    \[\phi\colon E(\Q_p)\times X(\Q_p)\longrightarrow X(\Q_p)\times X(\Q_p),\quad (T,x)\mapsto(T+x,x).\]
    For any $x\in X(\Q_p)_1$, the set
    \[V\colonequals\phi^{-1}(X(\Q_p)_1\times X(\Q_p)_1)=\brk{(T,x):x,T+x\in X(\Q_p)_1}\subset E(\Q_p)\times X(\Q_p)\]
    is an open neighbourhood of $(0,x)$, so there exists a small open subgroup $H_x\le E(\Q_p)$ and a small open set $U_x\subset X(\Q_p)_1$ such that $(0,x)\in H_x\times U_x\subset V$. As $x$ varies, the $U_x$'s cover the compact $X(\Q_p)_1$, so there is a finite subcollection of them -- say $x_1,\dots,x_n$ -- such that $X(\Q_p)_1=\bigcup_{i=1}^nU_{x_i}$. Hence, $\bigcap_{i=1}^nH_{x_i}$ is an open subgroup contained in $\Lambda_p$. This completes the verification that the subgroup $\Lambda_p\subset E(\Q_p)$ is open. Thus, $E(\Q_p)/\Lambda_p$ is discrete and compact, so $\Lambda_p$ has finite index.

    \emph{Part (ii): Global subgroups.} Consider the natural homomorphism
    \[\theta\colon E(\Q)\longrightarrow\prod_{p\mid\bad_F}E(\Q_p)/\Lambda_p.\]
    By the above, the codomain is finite, so $\ker\theta$ has finite index in $E(\Q)$. We claim that
    \begin{equation}\label{eq:local-coset-stable}
        X(\Q)_1+\ker\theta\subset X(\Q)_1.
    \end{equation}
    Let $P\in X(\Q)_1$ and $T\in\ker\theta$, and fix $p\mid\bad_F$. Then $P\in X(\Q_p)_1$ by \cref{eq:X1-local}, and $T$ maps to $\Lambda_p$ by construction, so $T+P\in X(\Q_p)_1$. As this holds for every $p\mid\bad_F$, we conclude that $T+P\in X(\Q)_1$ by \cref{eq:X1-local}. 

    \ref{it:Lambda1-finite-index} \Cref{eq:local-coset-stable} shows that $\ker\theta\le\Lambda_1$; since $\ker\theta$ is finite-index, so is $\Lambda_1$.

    \ref{it:Lambda1-cosets} By definition, $X(\Q)_1+\Lambda_1=X(\Q)_1$, so $X(\Q)_1$ is a disjoint union of orbits of $\Lambda_1$. Since $E(\Q)$ acts simply transitively on $X(\Q)$ and $\Lambda_1$ is finite index in $E(\Q)$, there can only be finitely many such orbits.
\end{proof}

\section{Proof of \texorpdfstring{\Cref{thm:main}}{Theorem 1.1}}
\label{sec:proof} We start with three auxiliary lemmas.

\begin{lemma}
  \label[lemma]{lemma:Omega-is-phi-X1}
  $\Omega_F = \phi\paren{X(\Q)_1}$.
\end{lemma}

\begin{proof}
  By definition, $\Omega_F$ is the set of points
  $(-c_1x_1^{e_1}:c_3x_3^{e_3}) \in \Pone(\Q)$ with
  $(x_1,x_2,x_3) \in \mcU(\Z)$. By \cref{eq:phi-belyi} this point is
  $\phi(\pr(x_1,x_2,x_3))$, so
  $\Omega_F = \phi\paren{\pr(\mcU(\Z))} = \phi\paren{X(\Q)_1}$ by
  \cref{eq:primitive-locus}.
\end{proof}

\begin{lemma}
  \label[lemma]{lemma:X1-aut-stable}
  Suppose that $\bfe$ is Euclidean. Then the primitive locus $X(\Q)_1$ is
  stable under $\Aut(\phi_{\Q})$.
\end{lemma}

\begin{proof}
  Every Euclidean signature has some $w_i = 1$
  (\Cref{table:euclidean-weights}), so by \Cref{lemma:aut-rational} every
  rational automorphism $\alpha$ of $\phi$ is a sign change
  $\x_j \mapsto \pm\x_j$. If $P = \pr(x_1,x_2,x_3)$ with
  $(x_1,x_2,x_3) \in \mcU(\Z)$, then $\alpha(P) = \pr(\pm x_1,\pm x_2,\pm x_3)$,
  and $(\pm x_1,\pm x_2,\pm x_3)$ is again a primitive solution. Hence
  $\alpha(P) \in X(\Q)_1$.
\end{proof}

Next we name the density contributed by the sieve.

\begin{definition}
  \label{def:rho}
  Suppose that $\bfe$ is Euclidean, let $E\colonequals\Jac(X)$, and let $\Lambda_1\subset E(\Q)$ be the primitive subgroup (\Cref{def:primitive-subgroup}).
  The \cdef{primitive density} of $F$ is
  \begin{equation}
    \label{eq:rho}
    \varrho_F \colonequals
    \frac{\#\brk{\text{orbits of } \Lambda_1 \text{ acting on } X(\Q)_1}}
    {[E(\Q):\Lambda_1]}.
  \end{equation}
\end{definition}

We restate our main theorem (\Cref{thm:main}) with the notation of
\Cref{sec:minimal-solutions} and \Cref{sec:proof}, and prove it.
\begin{theorem}
  \label{thm:main-precise}
  Suppose that $\bfe$ is Euclidean and that $\Omega_F \neq \emptyset$. 
  Let $E\coloneqq\Jac(X)$, and let $r = r(F)$ be the rank of $E(\Q)$.
  \begin{enumerate}[leftmargin=*, label=(\alph*)]
  \item If $r = 0$, then $\Omega_F$ is finite, and $N(\Omega_F;H) = \#\Omega_F$
    for all $H$ sufficiently large.
  \item If $r > 0$, then
    \begin{equation}
      \label{eq:main-asymptotic}
      N(\Omega_F;H) = \kappa(F)\cdot \paren{\log H}^{r/2}
      + O\paren{\paren{\log H}^{(r-1)/2}}, \qquad H \to \infty,
    \end{equation}
    where the constant is positive and given by
    \begin{equation}
      \label{eq:kappa}
      \kappa(F) =  \varrho_F\cdot \frac{\delta(\phi)}{\#\Aut(\phi_{\Q})}
      = \varrho_F\cdot\frac{\#E(\Q)_{\tors}}
      {\#\Aut(\phi_{\Q})\cdot\Gamma\paren{1+\tfrac r2}\sqrt{\Reg(E/\Q)}}
      \cdot \paren{\dfrac{2\pi\ell}{e_1e_2e_3}}^{r/2}.
    \end{equation}
  \end{enumerate}
  Here $\varrho_F$ is the primitive density of \Cref{def:rho}, and $\Reg(E/\Q)$ is the regulator
  with respect to $\OO_E(2O)$ (\Cref{remark:regulator-normalisation}).
\end{theorem}

\begin{proof}
  If $r = 0$ then $E(\Q)$ is finite, hence so is
  $\Omega_F = \phi\paren{X(\Q)_1}$, and (a) is immediate. Assume from now on
  that $r > 0$.

  Since $\Omega_F \neq \emptyset$, we may fix a point $O \in X(\Q)$. The map
  $P \mapsto [P - O]$ is an isomorphism $X \to E$ that is equivariant for
  translation by $E$. We use it to regard $\phi$ as a finite map
  $E \to \Pone_\Q$ of degree $e_1e_2e_3/\ell$ (\Cref{lemma:degree}), so that the
  results of \Cref{sec:counting} apply, and the orbits in $X(\Q)$ of a subgroup
  of $E(\Q)$ become its cosets in $E(\Q)$.

  Let $\Lambda_1 \subseteq E(\Q)$ be the primitive subgroup. By
  \Cref{thm:primitive-subgroup}\ref{it:Lambda1-cosets},
  \[
    X(\Q)_1 = \bigsqcup_{k=1}^{n} \paren{P_k + \Lambda_1}, \qquad
    n \colonequals \#\brk{\text{orbits of }\Lambda_1\text{ acting on } X(\Q)_1}.
  \]
  By \cref{eqn:delta-phi},
  $\delta(\phi;\Lambda_1) = \delta(\phi)/[E(\Q):\Lambda_1]$. Applying
  \Cref{prop:neron-corollary} with $\Lambda = \Lambda_1$ and $P_0 = P_k$ for
  each $k = 1, \dots, n$, and adding, we obtain
  \begin{equation}
    \label{eq:count-on-X1}
    \#\brk{P \in X(\Q)_1 \;:\; \Ht(\phi(P)) \leq H}
    = \frac{n\,\delta(\phi)}{[E(\Q):\Lambda_1]}\paren{\log H}^{r/2}
    + O\!\paren{\paren{\log H}^{(r-1)/2}},
  \end{equation}
  and the leading coefficient is $\varrho_F\,\delta(\phi)$ by \cref{eq:rho}.

  It remains to pass from $X(\Q)_1$ to its image
  $\Omega_F = \phi\paren{X(\Q)_1}$. The map $\phi$ is geometrically Galois
  (\Cref{prop:aut}), and $X(\Q)_1$ is stable under $\Aut(\phi_\Q)$ by \Cref{lemma:X1-aut-stable}, so \Cref{lem:counting-image-psi} with $S = X(\Q)_1$ implies that
  \[
    N(\Omega_F;H) = \frac{1}{\#\Aut(\phi_{\Q})}
    \paren{\#\brk{P \in X(\Q)_1 \;:\; \Ht(\phi(P)) \leq H} + O(1)}.
  \]
  Since $r \geq 1$, the $O(1)$ term is absorbed into the error term of
  \cref{eq:count-on-X1}, and we obtain \cref{eq:main-asymptotic} with
  $\kappa(F) = \varrho_F\,\delta(\phi)/\#\Aut(\phi_{\Q})$.

  The second expression for $\kappa(F)$ in \cref{eq:kappa} follows from
  \cref{eq:delta-explicit} and $\deg\phi = e_1e_2e_3/\ell$
  (\Cref{lemma:degree}). Finally, $\kappa(F) > 0$ because $\delta(\phi) > 0$ by
  \cref{eqn:delta-phi} and $\varrho_F > 0$ by definition.
\end{proof}

Combining \Cref{lemma:aut-rational,lemma:delta-explicit,thm:main-precise},
the constant reads as follows in the three Euclidean signatures.

\begin{table}[htbp]
  \centering
  \caption{Leading constant $\kappa(F)$ for Euclidean $F$ in
    \Cref{thm:main-precise} when $r(F) > 0$. 
    Here, $E$ denotes the Jacobian of $X$.
    For $\bfe = (3,3,3)$ one has $\varrho_F = 1$
    (\Cref{remark:333-no-sieve}).}
  \label{table:kappa}
  \rowcolors{2}{white}{gray!15}
  % \fracstrut gives each row room for the displayed fraction in the last column
  \newcommand{\fracstrut}{\rule[-3.2ex]{0pt}{8.4ex}}
  \begin{tabular}{cccc}
    \toprule
    $\bfe$ & $\deg \phi$ & $\#\Aut(\phi_{\Q})$ & $\kappa(F)$ \\
    \midrule
    $(3,3,3)$ & $9$ & $1$ & $
                            \dfrac{\#E(\Q)_{\tors}\cdot\paren{2\pi/9}^{r/2}}{\Gamma(1+\tfrac
                            r2)\cdot \sqrt{\Reg(E/\Q)}}$\fracstrut \\
    $(3,2,6)$ & $6$ & $2$ & $ \varrho_F  \cdot \dfrac{\#E(\Q)_{\tors}\cdot
                            \paren{\pi/3}^{r/2}}{2\cdot\Gamma(1+\tfrac r2)\cdot\sqrt{\Reg(E/\Q)}}$\fracstrut \\
    $(4,2,4)$ & $8$ & $4$ & $ \varrho_F \cdot\dfrac{\#E(\Q)_{\tors}\cdot
                            \paren{\pi/4}^{r/2}}{4\cdot\Gamma(1+\tfrac r2)\cdot\sqrt{\Reg(E/\Q)}}$\fracstrut \\
    \bottomrule
  \end{tabular}
\end{table}

\begin{remark}[The homogeneous case]
  \label{remark:333-no-sieve}
  For the signature $\bfe = (3,3,3)$ one has $\bfw = (1,1,1)$, so
  $\bad_F = 1$ and there are no bad primes. Hence $X(\Q)_1 = X(\Q)$ by
  \cref{eq:X1-local}, and $\varrho_F = 1$.
\end{remark}

\begin{remark}[The primitivity sieve can remove everything]
  \label{remark:sieve-can-kill}
  Let $F = \x_1^4 - 8\x_2^2 + \x_3^4$ be as in
  \Cref{ex:sameera's-example-1}. Then $(2:2:2) \in X(\Q)$, and yet
  $\Omega_F = \emptyset$. 
  
  Let $(x_1,x_2,x_3)\in\Z^3$ be a solution. Then, $x_1^4+x_3^4 \equiv 0\pmod8$, so
  $x_1,x_3$ must both be even. Hence, $16$ divides $8x_2^2=x_1^4+x_3^4$, so
  $x_2$ must be even as well, meaning $(x_1,x_2,x_3)$ cannot be primitive. In
  fact, this even shows $X(\Q_2)_1=\emptyset$ and so $\varrho_F=0$. Consequently, the
  hypothesis $\Omega_F\neq\emptyset$ in \Cref{thm:main} is \emph{strictly stronger} than
  ``$X$ has a rational point''.

\end{remark}

\begin{remark}[Normalization of the regulator]
  \label{remark:regulator-normalisation}
  In \cref{eq:delta-explicit,eq:kappa} the regulator $\Reg(E/\Q)$ is taken with
  respect to the N\'eron--Tate height attached to the line bundle $\OO(2O)$.
  Tables that use the height attached to $\OO(O)$ instead differ by a factor
  $2^r$ in the regulator.
\end{remark}

\bibliographystyle{amsalpha}
\bibliography{references.bib}
\end{document}